\documentclass[12pt,reqno,psfig,openbib]{amsart}

\usepackage[pagebackref=true, colorlinks=true, citecolor=blue]{hyperref}
\usepackage{amssymb,amsmath,graphicx,amsfonts,euscript}
\usepackage{fancyhdr}
\usepackage{epsfig}
\usepackage{setspace}
\usepackage{dsfont}
\usepackage[all]{xy}
\usepackage{amsmath}
\usepackage{enumitem}
\usepackage{algorithm}

\usepackage{caption}
\usepackage{subcaption}

\usepackage{soul}

\usepackage{float}

\usepackage{color}
\definecolor{darkgreen}{rgb}{0.0, 0.5, 0.0}
 \usepackage{pst-grad} 
 \usepackage{pst-plot} 

\numberwithin{equation}{section}
\numberwithin{figure}{section}

 \usepackage[left=1.05in,top=1.15in,right=1.05in]{geometry}
\usepackage{verbatim}

\newcommand{\spacer}{\vspace{1.5mm}}

\date{\today}

\renewcommand{\epsilon}{\varepsilon}
\newcommand{\ep}{\varepsilon}

\newcommand{\pa}{\partial}
\newcommand{\tht}{\theta}

\newcommand{\kap}{\kappa}
\newcommand{\blds}{\boldsymbol}


\begin{document}

\title[SL-PINN methods for Burgers' equation]
{Singular layer PINN methods for Burgers' equation}

\author[T.-Y. Chang, 
G.-M. Gie, 
Y. Hong, and 
C.-Y.  Jung]
{
Teng-Yuan Chang$^{1}$, 
Gung-Min Gie$^{2}$, 
Youngjoon Hong$^{1}$, and   
Chang-Yeol  Jung$^3$}
\address{$^1$ Department of Mathematical Sciences, KAIST, Korea}
\address{$^2$ Department of Mathematics, University of Louisville, Louisville, KY 40292}
\address{$^3$ Department of Mathematical Sciences, Ulsan National Institute of Science and Technology, 
Ulsan 44919, Korea}
\email{ty.chang@kaist.ac.kr}
\email{gungmin.gie@louisville.edu}
\email{hongyj@kaist.ac.kr}
\email{cjung@unist.ac.kr}

\begin{abstract}
In this article, we present a new learning method called sl-PINN to tackle the one-dimensional viscous Burgers problem at a small viscosity, which results in a singular interior layer. To address this issue, we first determine the corrector that characterizes the unique behavior of the viscous flow within the interior layers by means of asymptotic analysis. We then use these correctors to construct sl-PINN predictions for both stationary and moving interior layer cases of the viscous Burgers problem. Our numerical experiments demonstrate that sl-PINNs accurately predict the solution for low viscosity, notably reducing errors near the interior layer compared to traditional PINN methods. Our proposed method offers a comprehensive understanding of the behavior of the solution near the interior layer, aiding in capturing the robust part of the training solution.
\end{abstract}

\maketitle

\tableofcontents

\section{Introduction}
We consider the 1D Burgers' equation, particularly focusing on scenarios where the viscosity is small, resulting in the formation of interior layers,
\begin{equation}\label{Burgers}
    \left\{
    \begin{array}{rl}
        \spacer
        \pa_t u^\ep - \ep \pa_x^2 u^\ep + u^\ep \, \pa_x u^\ep = 0, 
            &  x \in \mathbb{R}, \,\, t > 0,\\
        u^\ep|_{t=0} = u_0, 
            & x \in \mathbb{R},
    \end{array}
    \right.
\end{equation}
where $u_0 = u_0(x)$ is the given initial data and
$\ep>0$ is the small viscosity parameter. 
The Burgers' equation (\ref{Burgers}) is 
supplemented with one of the following boundary conditions, 
\begin{itemize}
    \item[($i$)] 
            $u^\ep$ is periodic in $-1 < x < 1$,  
            provided that the initial data $u_0$ in is periodic in $(-1, 1)$ as well;
    \item[($ii$)]
            $u^\ep \rightarrow g_L(t)$ as $x \rightarrow - \infty$ and  
            $u^\ep \rightarrow g_R(t)$ as $x \rightarrow  \infty$ where 
            the boundary data $g_L$ and $g_R$ are smooth functions in time.         
\end{itemize}
In the examples presented in Section \ref{s.num}, we consider the case where the data $u_0$, $g_L$, and $g_R$ are sufficiently regular. Therefore, the Burgers' equation at a fixed viscosity $\epsilon > 0$ is well-posed up to any finite time $T > 0$ and satisfies  
\begin{equation}\label{e:bound_viscous_sol}
    |u^\ep (x, t)| 
        \leq C_0, 
        \qquad 
        \forall (x, t) \in \mathbb{R}\times (0, T),
\end{equation}
for a constant $C_0> 0 $ depending on the data, but independent of the viscosity parameter $\ep>0$; see, e.g., \cite{CJL23} for the proof when the boundary condition ($ii$) above is imposed.  

The corresponding limit $u^0$ of $u^\ep$ at the vanishing viscosity $\ep = 0$ satisfies the inviscid Burgers' equation,  
\begin{equation}\label{Burgers_inviscid}
    \left\{
    \begin{array}{rl}
        \spacer
        \pa_t u^0 + u^0 \, \pa_x u^0 = 0, 
            &  x \in \mathbb{R}, \,\, t > 0,\\
        u^0(x,0) = u_0 (x), 
            & x \in \mathbb{R},
    \end{array}
    \right.
\end{equation}
supplemented with ($i$) or ($ii$) above as for the viscous problem (\ref{Burgers}).  
Concerning all the examples in Section \ref{s.num}, 
the limit solution $u^0$ exhibits a discontinuous interior shock  
located at $x = \alpha(t)$, $t\geq0
$ (either stationary $d \alpha / dt = 0$ or moving $d \alpha / dt \neq 0$), and $u^0$ is smooth on each side of $x=\alpha(t)$.  
Hence by identifying the states of a shock, we introduce smooth one-side solutions, 
\begin{equation}\label{e:u^0_LR}
\begin{array}{l}
\spacer
    u^- 
        = u^-(t)
        = \lim_{x \rightarrow \alpha(t)^-}u^0 (x, t), \\
    u^+ 
        = u^+(t)
        = \lim_{x \rightarrow \alpha(t)^+}u^0 (x, t).
\end{array}
\end{equation}
Setting $\dot{\alpha}(t) = d \alpha / dt$, 
the solution $u^0$ satisfies that 
\begin{equation}\label{e:COND}
    u^- > \dot{\alpha}(t) > u^+, 
    \quad
    t \geq 0,
\end{equation}
and, in addition, we have 
\begin{equation}\label{e:COND_alpha}
    \dot{\alpha}(t) 
    =
        \dfrac{u^- + u^+}{2}, 
    \quad
    t \geq 0.
\end{equation}
Note that, in this case, the viscous solution $u^\ep$ is located between the left and right solutions of $u^0$ at the shock $x = \alpha(t)$, i.e., 
\begin{equation}\label{e:u^epVSu^0}
    u^-
    >
    \dfrac{u^- + u^\ep}{2}
    >
    \dot{\alpha}(t) 
    =
        \dfrac{u^- + u^+}{2}
    >
    \dfrac{u^\ep + u^+}{2}
    >
    u^+, 
    \quad
    \text{at }
    x = \alpha(t). 
\end{equation}



The primary goal of this article is to analyze and approximate the solution $u^\ep$ of (\ref{Burgers}) as the viscosity $\ep$ tends to zero 
when the (at least piecewise) smooth initial data creates an interior shock of $u^0$.  
Specifically, we introduce a new machine learning technique called a {\em singular-layer Physics-Informed Neural Network} (sl-PINN) to accurately predict solutions to (\ref{Burgers}) when the viscosity is small.
Following the methodology outlined in \cite{GHJ24, GHJM24, GHJL24, CGHJ24}, in Section \ref{S.Asymp}, we start by revisiting the asymptotic analysis for (\ref{Burgers}) as presented in \cite{CJL23}; see, e.g., \cite{CJL19, GJL22, Book16, SK87} as well. Then, utilizing the interior layer analysis, we proceed to build our new sl-PINNs for (\ref{Burgers}) in Sections \ref{PINN_Burgers}. The numerical computations for our sl-PINNs, along with a comparison to those from conventional PINNs, are detailed in Section \ref{s.num}. Finally, in the conclusion in Section \ref{S.conclusion}, we summarize that the predicted solutions via sl-PINNs approximate the stiff interior layer of the Burgers problem quite accurately.

\section{Interior layer analysis}\label{S.Asymp}
We briefly recall from \cite{CJL23} the interior layer analysis for (\ref{Burgers}) that is made useful in the section below \ref{PINN_Burgers} where we construct our {\em singular-layer Physics-Informed Neural Networks} ({\em sl-PINNs}) for Burgers' equation.
To analyze the interior layer of the viscous Burgers' equation (\ref{Burgers}), 
we first introduce a moving coordinate in space, centered at the shock location $\alpha(t)$,  
\begin{equation}\label{e:moving_coordinate}
    \tilde{x} 
        = \tilde{x}(t)
        = x - \alpha(t), 
        \quad
        t \geq 0, 
\end{equation}
and write the solutions of (\ref{Burgers}) and (\ref{Burgers_inviscid}) with respect to the new system $(\tilde{x}, t)$ in the form, 
\begin{equation}\label{e:sols_new_coordinate}
    u^\ep(x, t)
        =   \tilde{u}^\ep (\tilde{x}, t), 
    \qquad
    u^0(x, t)
        =   \tilde{u}^0 (\tilde{x}, t).  
\end{equation}
The initial condition is written in terms of the new coordinate as   
\begin{equation}\label{e:initial_new_coordinate}
    u_0(x, t)
        =   \tilde{u}_0 (\tilde{x}, t).  
\end{equation}
Now we set 
\begin{equation}\label{e:u^ep_decompose}
u^\ep(x, t)
    =
    \tilde{u}^\ep (\tilde{x}, t) 
    = 
    u^{\ep}_{L}(\tilde{x}, t)  + u^{\ep}_{R}(\tilde{x}, t) , 
\end{equation}
where $u^{\ep}_{L}$ and $u^{\ep}_{R}$ are the 
disjoint and smooth portions of $u^\ep$ 
located on the left-hand side and the right-hand side of the shock at  $x=\alpha(t)$, 
and they are 
defined by
\begin{equation}\label{e:u_decompose_2}
\left\{
\begin{array}{l}
\spacer
    u^\ep_L (\tilde{x}, t)
        =   
            \tilde{u}^\ep (\tilde{x}, t), 
            \quad 
            \text{for } -\infty < \tilde{x} < 0, 
            \, 
            t \geq 0,\\
    u^\ep_R (\tilde{x}, t)
        =   
            \tilde{u}^\ep (\tilde{x}, t), 
            \quad 
            \text{for } 0< \tilde{x} < \infty, 
            \, 
            t \geq 0.
\end{array}
\right.
\end{equation}

Using (\ref{Burgers}) and (\ref{e:moving_coordinate}), and applying the chain rule, 
one can verify that the $u^\ep_L$ and $u^\ep_R$ satisfy the equations below:
\begin{equation}\label{BurgersL}
    \left\{
    \begin{array}{rl}
        \spacer
        \pa_t u^\ep_L 
        - \ep \pa_{\tilde{x}}^2 u^\ep_L 
        + 
        \pa_{\tilde{x}} \left( \left( \displaystyle \frac{u^\ep_L}{2} - \dot{\alpha}(t) \right) u^\ep_L \right) = 0, 
            &  -\infty < \tilde{x} < 0, \,\, t > 0,\\
        \spacer
        u^\ep_L(0, t) 
            = \tilde{u}^\ep(0, t) 
            = u^\ep(\alpha(t), t), 
            & \tilde{x} = 0,\\
        \spacer
        u^\ep_L \to g_L(t), 
            & \text{as } \tilde{x} \to -\infty,\\
        u^\ep_L(\tilde{x},0) 
        = \tilde{u}_0 (\tilde{x}(0))
        = u_0 (\tilde{x} + \alpha(0))
        , 
            & -\infty < \tilde{x} < 0,
    \end{array}
    \right.
\end{equation}
and 
\begin{equation}\label{BurgersR}
    \left\{
    \begin{array}{rl}
        \spacer
        \pa_t u^\ep_R 
        - \ep \pa_{\tilde{x}}^2 u^\ep_R 
        + \pa_{\tilde{x}} \left( \left( \displaystyle \frac{u^\ep_R}{2} - \dot{\alpha}(t) \right) u^\ep_R \right) = 0, 
            &  0 < \tilde{x} < \infty, \,\, t > 0,\\
        \spacer
        u^\ep_R(0, t) 
        = \tilde{u}^\ep(0, t), 
        = u^\ep(\alpha(t), t), 
            & \tilde{x} = 0,\\
        \spacer
        u^\ep_R \to g_R(t), 
            & \text{as } \tilde{x} \to \infty,\\
        u^\ep_R(\tilde{x},0) 
        = \tilde{u}_0 (\tilde{x}(0))
        = u_0(\tilde{x} + \alpha(0)),
            & 0 < \tilde{x} < \infty.
    \end{array}
    \right.
\end{equation}
\smallskip

Thanks to the moving coordinate system $(\tilde{x}, t)$, 
now we are ready to analyze the interior shock layer near $x = \alpha(t)$ by investigating the two disjoint boundary layers 
of $u^\ep_L$ and $u^\ep_R$,  
located on the left-hand and right-hand sides of $x = \alpha(t)$. 
To this end, we introduce asymptotic expansions in the form, 
\begin{equation}\label{e:expansions}
\left\{
    \begin{array}{l}
         \spacer
            u^\ep_L(\tilde{x}, t)
                =
                    u^0_L(\tilde{x}, t) + \varphi_L(\tilde{x}, t), 
            \quad 
            -\infty < \tilde{x} < 0, \,\, t > 0,\\
            u^\ep_R(\tilde{x}, t)
                =
                    u^0_R(\tilde{x}, t) + \varphi_R(\tilde{x}, t), 
            \quad
            0 < \tilde{x} < \infty, \,\, t > 0.
    \end{array}
\right.
\end{equation}
Here 
$u^{0}_{L}$ and $u^{0}_{R}$ are the 
disjoint and smooth portions of $u^0$ 
located on the left-hand side and the right-hand side of $x=\alpha(t)$,  
given by 
\begin{equation}\label{e:u^0_decompose_2}
\left\{
\begin{array}{l}
\spacer
    u^0_L (\tilde{x}, t)
        =   
            \tilde{u}^0 (\tilde{x}, t), 
            \quad 
            \text{for } -\infty < \tilde{x} < 0, 
            \, 
            t \geq 0,\\
    u^0_R (\tilde{x}, t)
        =   
            \tilde{u}^0 (\tilde{x}, t), 
            \quad 
            \text{for } 0< \tilde{x} < \infty, 
            \, 
            t \geq 0, 
\end{array}
\right.
\end{equation}
i.e., $\tilde{u}^0 = u^0_L + u^0_R$ for each $(\tilde{x}, t) \in \mathbb{R} \times (0, \infty)$.

The $\varphi_L$ or $\varphi_R$ in the asymptotic expansion  (\ref{e:expansions}) is an artificial function, called the  corrector, 
that approximates the stiff behavior of $u^\ep_* - u^0_*$, $* = L, R$, near the shock of $u^0$ at $\tilde{x} = 0$ (or equivalently $x = \alpha(t)$). 
Following the asymptotic analysis in \cite{CJL23}, 
we define the correctors $\varphi_L$ and  $\varphi_R$ as the solutions of   
\begin{equation}\label{corL}
    \left\{
    \begin{array}{rl}
        \spacer \displaystyle
        -\pa_{\tilde{x}}^2 {\varphi}_L 
        + \frac{1}{2} \pa_{\tilde{x}} \left( 2 u^0_L(0,t) {\varphi}_L 
        - 2 \dot{\alpha}(t) {\varphi}_L 
        + ({\varphi}_L)^2 \right) = 0, 
            &  -\infty < \tilde{x} < 0, \,\, t > 0,\\
        \spacer
        {\varphi}_L(0, t) = b(t) - u^0_L(0, t), 
            & \tilde{x} = 0,\\
        {\varphi}_L \to 0, 
            & \text{as } \tilde{x} \to -\infty,
    \end{array}
    \right.
\end{equation}
and
\begin{equation}\label{corR}
    \left\{
    \begin{array}{rl}
        \spacer \displaystyle
        - \ep \pa_{\tilde{x}}^2 {\varphi}_R 
        + \frac{1}{2} \pa_{\tilde{x}} \left( 2 u^0_R(0,t) {\varphi}_R 
        - 2 \dot{\alpha}(t) {\varphi}_R + ({\varphi}_R)^2 \right) = 0, 
            &  0 < \tilde{x} < \infty, \,\, t > 0,\\
        \spacer
        {\varphi}_R(0, t) = b(t) - u^0_R(0, t), 
            & \tilde{x} = 0,\\
        {\varphi}_R \to 0, 
            & \text{as } \tilde{x} \to \infty,
    \end{array}
    \right.
\end{equation}
where 
\begin{equation}\label{e:b}
    b 
        = b(t)
        = u^{\ep}(\alpha(t), t).
\end{equation}
Thanks to the property of $u^\ep$, (\ref{e:bound_viscous_sol}), 
we notice that 
the value of $b$ is bounded by a constant 
independent of the small viscosity parameter $\ep$.

We recall from Lemma 3.1 in \cite{CJL23} that the explicit expression of $\varphi_L$ (and hence of $\varphi_R$ by the symmetry) is given by  
    \begin{align}\label{e:corrector_explicit}
        {\varphi}_* \left( {\tilde{x}}, t \right)
        &= \left( \dot{\alpha}(t) - u^0_*(0,t) \right) \left( 1 +  \tanh \left( \frac{M_*}{2} (\tilde{x}-m_*) \right) \right),
    \quad 
    * = L, R, 
    \end{align}
    where
    \begin{align*}
        M_* 
            = 
            \dfrac{u^0_*(0,t)-\dot{\alpha}(t)}{\ep}, 
        \qquad  
        m_*
            =
            \dfrac{1}{M_*}\log \left( \frac{u^0_*(0,t)-2\dot{\alpha}(t)+b(t)}{u^0_*(0,t)-b(t)} \right). 
    \end{align*}


The validity of our asymptotic expansion of $u^\ep$ in (\ref{e:u^ep_decompose}) and (\ref{e:expansions}) as well as the vanishing viscosity limit of $u^\ep$ to $u^0$ 
are rigorously verified in \cite{CJL23}, and here 
we briefly recall the result without providing the detailed and technical proof; see Theorem 3.5 in \cite{CJL23}:

Under a certain regularity condition on the limit solution $u^0_L$ and $u^0_R$, the asymptotic expansion in (\ref{e:u^ep_decompose}) and (\ref{e:expansions}) is valid in the sense that 
\begin{equation}\label{e:validity_expansion}
\left\{
\begin{array}{l}
\spacer
    \|
        u^\ep_L - (u^0_L + \varphi_L)
    \|_{L^\infty(0, T; L^2(\mathbb{R}^-))}
    +
    \ep^{\frac{1}{2}}
    \|
        u^\ep_L - (u^0_L + \varphi_L)
    \|_{L^2(0, T; H^1(\mathbb{R}^-))}
        \leq 
            \kap \ep, 
    \\
    \|
        u^\ep_R - (u^0_R + \varphi_R)
    \|_{L^\infty(0, T; L^2(\mathbb{R}^+))}
    +
    \ep^{\frac{1}{2}}
    \|
        u^\ep_R - (u^0_R + \varphi_R)
    \|_{L^2(0, T; H^1(\mathbb{R}^+))}
        \leq 
            \kap \ep,
\end{array}
\right.
\end{equation}
for a generic constant $\kap>0$ depending on the data, but independent of the viscosity parameter $\ep$. 
Moreover, the viscous solution $u^\ep$ converges to the corresponding limit solution $u^0$ as the viscosity vanishes in the sense that 
\begin{equation}\label{e:VVL}
    \left\{
    \begin{array}{l}
        \spacer
        \|
        u^\ep_L - u^0_L 
        \|_{L^\infty(0, T; L^2(\mathbb{R}^-))}
        +
        \|
        u^\ep_R - u^0_R 
        \|_{L^\infty(0, T; L^2(\mathbb{R}^+))}
            \leq \kap \ep^{\frac{1}{2}},\\
        \|
        u^\ep_L - u^0_L 
        \|_{L^\infty(0, T; L^1(0, K))}
        +
        \|
        u^\ep_R - u^0_R 
        \|_{L^\infty(0, T; L^1(-K, 0)}
            \leq \kap \ep, 
        \quad \forall K>0, 
    \end{array}
    \right.
\end{equation}
for a  constant $\kap>0$ depending on the data, but independent of $\ep$.

\section{{PINNs and sl-PINNs  
for Burgers' equation
}}\label{PINN_Burgers}

In this section, we are developing a new machine learning process called the {\em singular-layer Physics-Informed Neural Network} (sl-PINN) for the Burgers' equation (\ref{Burgers}). 
We will be comparing its performance with that of the usual PINNs, specifically when the viscosity $\ep>0$ is small and hence  the viscous Burgers' solution shows rapid changes in a small region known as the interior layer.  
To start, we provide a brief overview of Physics-Informed Neural Networks (PINNs) for solving the Burgers equation. 
We then introduce our new sl-PINNs and discuss their capabilities in comparison to the traditional PINNs. 
We see below in Section \ref{s.num} 
that our new sl-PINNs yield more accurate predictions, particularly when the viscosity parameter is small.


\subsection{{PINN structure}}\label{PINN}
We introduce an $L$-layer feed-forward neural network (NN) defined by the recursive expression: 
$\mathcal{N}^l$, $0 \leq l \leq L$ such that 
\begin{align*}
    \text{input layer: }& \mathcal{N}^{0}(\blds{x}) = \blds{x} \in \mathbb{R}^{N_{0}}, \\
    \text{hidden layers: }& \mathcal{N}^{l}(\blds{x}) = \sigma(\blds{W}^{l} \mathcal{N}^{l-1}(\blds{x}) + \blds{b}^{l}) \in \mathbb{R}^{N_{l}}, \quad 1\leq l \leq L-1, \\
    \text{output layer: }& \mathcal{N}^{L}(\blds{x}) = \blds{W}^{L} \mathcal{N}^{L-1}(\blds{x}) + \blds{b}^{L} \in \mathbb{R}^{N_{L}},
\end{align*}
where $\blds{x}$ is an input vector in $\mathbb{R}^{N_{0}}$, and $\blds{W}^{l} \in \mathbb{R}^{N_{l} \times N_{l-1}}$ and $\blds{b}^{l} \in \mathbb{R}^{N_{l}}$, $l=1,\ldots,L$, are the (unknown) weights and biases respectively. Here, $N^{l}$, $l=0,\ldots,L$, is the number of neurons in the $l$-th layer, and $\sigma$ is the activation function. 
We set $\blds{\tht}=\{ (\blds{W}^{l}, \blds{b}^{l}) \}_{l=1}^{L}$ as the collection of all the parameters, weights and biases. 

To solve the Burgers equation (\ref{Burgers}), 
we use the NN structure above with 
the output $ \mathcal{N}^{L} = \hat{u} (x,t;\blds{\tht})$ 
and 
the input $\blds{x} = (x, t) \in (-1, 1) \times (0, \infty)$ 
where  
a 
proper treatment 
is implemented for the boundary condition $(i)$ or $(ii)$ in  (\ref{Burgers}); see {(\ref{e:l_PINN_smooth}) and (\ref{e:l_PINN_nonsmooth})} below. 
In order to 
determine the predicted solution $\hat{u}$ of (\ref{Burgers}),  
we use the 
loss function defined by
\begin{equation}\label{e:loss_PINN}
    \begin{split}
        \mathcal{L}(\blds{\tht}; \mathcal{T}, \mathcal{T}_{b}, \mathcal{T}_{0}) 
        &= \frac{1}{|\mathcal{T}|} \sum_{(x,t)\in \mathcal{T}} \left| \pa_{t} \hat{u} - \ep \pa_{x}^{2} \hat{u} + \hat{u}\, \pa_{x} \hat{u}   \right|^{2} \\
        &+ \frac{1}{|\mathcal{T}_{b}|} {  
        \sum_{t \in \mathcal{T}_{b}}  l( \hat{u} )  
        }
        + \frac{1}{|\mathcal{T}_{0}|} \sum_{x \in \mathcal{T}_{0}} \left| \hat{u}|_{t=0}-u_{0} \right|^{2} ,
    \end{split}
\end{equation}
where $\mathcal{T} \subset \Omega_{x} \times \Omega_{t}$ is a set of training data in the space-time domain, $\mathcal{T}_{b} \subset \pa \Omega_{x} \times \Omega_{t}$ on the boundary, and $\mathcal{T}_{0} \subset \Omega_{x}$ on the initial boundary at $t=0$. 
The function $l(\hat{u})$ in (\ref{e:loss_PINN}) is defined differently for different examples below to enforce the boundary conditions $(i)$ or $(ii)$ in (\ref{Burgers}); see {(\ref{e:l_PINN_smooth}) and (\ref{e:l_PINN_nonsmooth})} below.

\subsection{{sl-PINN structure}}\label{sl-PINN}

In order to create our predicted solution for the Burgers equation (\ref{Burgers}), 
we revisit the interior layer analysis in Section \ref{S.Asymp} and then, using this analytic result, we develop two training solutions in distinct regions. 
More precisely, 
we partition the domain $\Omega_{x}$ into disjoint subdomains designated as $\Omega_{x,L}$ and $\Omega_{x,R}$, which are separated by the shock curve $x=\alpha(t)$. 
After that, we introduce the output functions $\hat{u}_{*}=\hat{u}_{*}(x,t; \blds{\tht}_{*})$ for the neural network in the subdomains $\Omega_{x,*}$ where $*=L,R$. The training solution for our sl-PINN method is then defined as:

\begin{equation} \label{tilde u}
    \tilde{u}_* (x,t; \blds{\tht}_*) 
        = 
            \hat{u}_* (x,t; \blds{\tht}_*) 
            + 
            \tilde{\varphi}_* (x,t; \blds{\tht}_*), 
    \quad
    *=L,R
    .
\end{equation}
Here $\tilde{\varphi}_{*}$ is the training form of the corrector $\varphi$ in (\ref{e:corrector_explicit}), 
and it is given by
\begin{equation} \label{cor}
    \tilde{\varphi}_*(x,t; \blds{\tht}_*) 
        = 
        (\dot{\alpha}(t)
        -
        \hat{u}_* (\alpha(t),t; \blds{\tht}_*) ) 
        \left( 
        1+\tanh\left( 
            \frac{\tilde{M_*}}{2} (x-\alpha(t) - \tilde{m} ) 
            \right) 
        \right),
\end{equation}
with the trainable parameters, 
\begin{equation}\label{e:paramemter_cor_NN}
    \tilde{M}_* 
        = 
        \frac{\hat{u}_* (\alpha(t),t; \blds{\tht}_*)
            -\dot{\alpha}(t)
        }{\ep}, 
    \quad 
    \tilde{m}_* 
        = 
        \frac{1}{\tilde{M}_*}
        \log \left( 
                \frac{
                \hat{u}_*(\alpha(t),t; \blds{\tht}_*)-2\dot{\alpha}(t)+\tilde{b}
                }
                {\hat{u}_*(\alpha(t),t; \blds{\tht}_*)-\tilde{b}
                } 
            \right). 
\end{equation}
Note that the $\tilde{b}$ appears in (\ref{e:paramemter_cor_NN}) is an approximations of $b(t) = u^{\ep}(\alpha(t),t)$ and 
{it is  updated iteratively during the network training}; see 
Algorithm \ref{algo} below. 

As noted in Section \ref{S.Asymp}, 
the correctors $\tilde{\varphi}_{*}$, $*= L, R$, 
capture  well the 
singular behavior of the viscous Burgers' solution $u^\ep$ near the shock of $u^0$ at $x = \alpha(t)$. 
{ Hence, by integrating these correctors into the training solution, we record the stiff behavior of $u^\ep$ directly in (\ref{tilde u}). 
}

To define the loss function for our sl-PINN method, 
we modify and use the idea introduced in the eXtended Physics-Informed Neural Networks (XPINNs), e.g., \cite{JK20}. 
In fact, to manage the mismatched value across the shock curve as well as the residual of differential equations along the interface, we define the loss function as   
\begin{equation} \label{loss slPINN}
    \begin{split}
        &\mathcal{L}_* (\tht_*; \mathcal{T}_{*}, \mathcal{T}_{*,b}, \mathcal{T}_{*,0}, \mathcal{T}_{\Gamma}) \\
        &= \frac{1}{|\mathcal{T}_{*}|}\sum_{(x,t) \in \mathcal{T}_{*}} \left| \pa_t \tilde{u}_* - \ep \pa_x^2 \tilde{u}_* + \tilde{u}_* \pa_x \tilde{u}_*  
        \right|^2 \\
        &+ 
        { 
        \frac{1}{|\mathcal{T}_{*,b}|}\sum_{t \in   
        \mathcal{T}_{*,b}}  l( \tilde{u}_* )  
        }
        + \frac{1}{|\mathcal{T}_{*,0}|} \sum_{x \in \mathcal{T}_{*,0}} \left| \tilde{u}_*|_{t=0} - u_{0} \right|^2 \\
        &+ \frac{1}{|\mathcal{T}_{\Gamma}|} \sum_{(x,t) \in \mathcal{T}_{\Gamma}} 
        {  
        | \tilde{u}_L - \tilde{u}_R |^2 
        }\\
        &+ \frac{1}{|\mathcal{T}_{\Gamma}|} \sum_{(x,t) \in \mathcal{T}_{\Gamma}} \left| 
        { 
        \left( 
        \pa_t \tilde{u}_L - \ep \pa_x^2 \tilde{u}_L + \tilde{u}_L \pa_x \tilde{u}_L  
        \right)  
        - 
        \left( 
        \pa_t \tilde{u}_R - \ep \pa_x^2 \tilde{u}_R + \tilde{u}_R \pa_x \tilde{u}_R
        \right) 
        }
        \right|^2, 
    \quad
        * =L, R.
    \end{split}
\end{equation}
{ 
Here the interface between $\Omega_L = (-\infty, \alpha(t))$ and $\Omega_R = (\alpha(t), \infty)$ is given by 
$\Gamma := \{(x=\alpha(t), \, t) \in \Omega_x \times \Omega_t \, | \,\, t \in \Omega_t \}$, 
}
and the training data sets are chosen such that $\mathcal{T}_{*} \subset \Omega_{x,*} \times \Omega_{t}$, $\mathcal{T}_{*,b} \subset \pa \Omega_{x,*} \times \Omega_{t}$, $\mathcal{T}_{*,0} \subset \Omega_{x,*} $, and $\mathcal{T}_{\Gamma} \subset \Gamma$. 
The function $l(\tilde{u}_*)$, $* = L, R$, in (\ref{loss slPINN}) is defined differently for different examples below to enforce the boundary conditions $(i)$ or $(ii)$ in (\ref{Burgers}); see {(\ref{e:l_slPINN_smooth}) and (\ref{e:l_slPINN_nonsmooth})} below.

The solutions $\tilde{u}_{L}$ and $\tilde{u}_{R}$ defined in (\ref{tilde u}) are trained simultaneously at each iteration step based on the corresponding loss functions $\mathcal{L}_{L}$ and $\mathcal{L}_{R}$ defined in (\ref{loss slPINN}). 
Note that we define two networks separately for $\tilde{u}_{L}$ and $\tilde{u}_{R}$, but 
they interact with each other via the loss function (\ref{loss slPINN}) at each iteration and thus 
the trainings of $\tilde{u}_{L}$ and $\tilde{u}_{R}$ occur simultaneously and they depend on each other.  
The training algorithm is listed below in Algorithm \ref{algo}.

\begin{algorithm} 
  \caption{(sl-PINN training algorithm) 
        }
\begin{flushleft}  
\begin{itemize}
    \item Construct two NNs on the subdomains $\Omega_{x,*} \times \Omega_{t}$, $\hat{u}_{*} = \hat{u}_{*}(x,t; \blds{\tht}^{(0)}_{*})$ for $*=L,R$, where $\blds{\tht}^{(0)}_{*}$ denotes the initial parameters of the NN.

    \item Initialize the $\tilde{b}$ by $
        \tilde{b} = 
        2^{-1} 
        \big( \lim_{x \rightarrow \alpha(0)^{-}}u_{0}(x) 
        +
        \lim_{x \rightarrow \alpha(0)^{+}}u_{0}(x) \big)$.

    \item Define the training solutions $\tilde{u}_{*} = \tilde{u}_{*}(x,t; \blds{\tht}^{(0)}_{*})$, $*=L,R$, as in (\ref{tilde u}).
    \item Choose the training sets of random data points as $\mathcal{T}_{*}$, $\mathcal{T}_{*,b}$, $\mathcal{T}_{*,0}$, and $\mathcal{T}_{\Gamma}$, $*=L,R$.
    \item Define the loss functions $\mathcal{L}_{*}$, $*=L,R$, as in (\ref{loss slPINN}).

    \item \textbf{for $k= 1,\ldots, \ell$  do}
    \begin{itemize}
        \item[]
            train $\tilde{u}_{*} =
            \tilde{u}_{*}^{(k)} =
            \tilde{u}_{*}(x,t; \blds{\tht}^{(k)}_{*})$ and obtain the updated parameters $\blds{\tht}_{*}^{(k+1)}$ for $* = L, R$

        \item[] 
            update $\tilde{b}$ by 
            $\tilde{b}
                =
                \frac{1}{2} \left(\tilde{u}_{L}\left(\alpha(t),t;\blds{\tht}_{L}^{(k+1)}\right)
                +
                \tilde{u}_{R}\left(\alpha(t),t;\blds{\tht}_{R}^{(k+1)}\right) \right)$

    \end{itemize}
    \item[] \textbf{end for}
    \item  
            Obtain the predicted solution as
            { 
            $$
            \tilde{u} 
                = 
                \tilde{u}^{(\ell)}_L + \tilde{u}^{(\ell)}_R.
            $$
            }
\end{itemize}
Here the superscript $(k)$ denotes the $k$-th iteration step by the optimizer, and 
{
$\ell$ is the number of the maximum iterations, or 
the one where an early stopping within a certain tolerance is achieved.}
\end{flushleft}
\label{algo}
\end{algorithm}

\section{Numerical Experiments} \label{s.num}

We are evaluating the performance of our new sl-PINNs in solving the viscous Burgers problem (\ref{Burgers}) in cases where the viscosity is small, resulting in a stiff interior layer. We are investigating this problem with both smooth initial data and non-smooth initial data. In the case of smooth initial data, we are using a sine function, which is a well-known benchmark experiment for the slightly viscous Burgers' equation \cite{CMP86, HM98, CP82}. Additionally, we are also considering initial data composed of Heaviside functions, which lead to the occurrence of steady or moving shocks (refer to Chapters 2.3 and 8.4 in \cite{Peter2014} or \cite{DL05, BZ18}).


We define 
the error between the 
exact solution and the predicted solution as
\begin{align*}
    E(x,t) = u_{ref}(x,t) - u_{pred}(x,t),
\end{align*}
where $u_{ref}$ is the high-resolution reference solution obtained by the finite difference methods, and $u_{pred}$ is the predicted solution obtained by either PINNs or sl-PINNs. 
The forward-time finite (central) difference method is employed to generate the reference solutions in all tests 
where 
the mesh size $\Delta x = 10^{-5}$ and the time step $\Delta t = 10^{-7}$ are chosen to be sufficiently small. 

When computing the error for each example below, we use the high-resolution reference solution instead of the exact solution. This is because the explicit expression, (\ref{e:exact_smooth}) or (\ref{e:exact_nonsmooth}), of the solution for each example is not suitable for direct use, as it is involved in integrations.

We measure below the error of each computation in the various norms,    
\begin{align}
    \| E \|_{L^{2}_{x}(\Omega)} &= \left( \int_{\Omega} \left|E(x,t)\right|^{2} dx \right)^{\frac{1}{2}}, \\
    \| E \|_{L^{2}_{t}( L^{2}_{x}(\Omega))} &= \left( \int_{0}^{T} \| E \|^{2}_{L^{2}_{x}(\Omega)} dt \right)^{\frac{1}{2}},  \label{L2L2}
\end{align}
and 
\begin{align} \label{Linf}
    \| E \|_{L^{\infty}_{x}(\Omega)} = \max_{x \in \Omega} \left|E(x,t)\right|. 
\end{align}

\subsection{Smooth initial data}\label{S.smooth_data} 
We consider the Burgers' equation (\ref{Burgers}) 
supplemented with with the periodic boundary condition $(i)$ 
where the initial data is given by 
\begin{align*}
    u_{0}(x) = -\sin (\pi x).
\end{align*}

For this example, the explicit expression of solution is given by 
\begin{align}\label{e:exact_smooth}
    u^\ep(x,t) = \dfrac{-\int_{\mathbb{R}} \sin(\pi(x-\eta)) \ F(x-\eta)\ G(\eta,t) \ d\eta}
    {\int_{\mathbb{R}} \ F(x-\eta)\ G(\eta,t) \ d\eta},
\end{align}
where
\begin{align*}
    F(y) = \exp \bigg(
            -\dfrac{\cos(\pi y)}{2\pi \ep}\bigg),  \quad 
    G(y,t) = \exp \bigg(
                    -\frac{y^{2}}{4\ep t}
                  \bigg).
\end{align*}
We have observed that even though the analytic solution above has been well-studied (refer to, for example, \cite{CMP86}), 
any numerical approximation of (\ref{e:exact_smooth}) requires a specific numerical method to compute the integration, and the integrand functions become increasingly singular when the viscosity $\ep$ is small. To avoid this technical difficulties for our simulations, 
we use a high-resolution finite difference method to produce our reference solution.


Now, to build our sl-PINN predicted solution, we first notice that 
the interior layer for this example is stationary, i.e., $\alpha(t)\equiv 0$, and hence we find that  
the training corrector (\ref{cor}) is reduced to
\begin{align*}
    \tilde{\varphi}_*(x,t; \blds{\tht}_*) 
        = 
            -\hat{u}_* (0,t; \blds{\tht}_*)  
            \left( 
            1+\tanh\left( \frac{\hat{u}_* (0,t; \blds{\tht}_*)}{2\ep} \ x \right) 
            \right).
\end{align*}

\begin{figure}
    \centering
    \includegraphics[width=9cm]{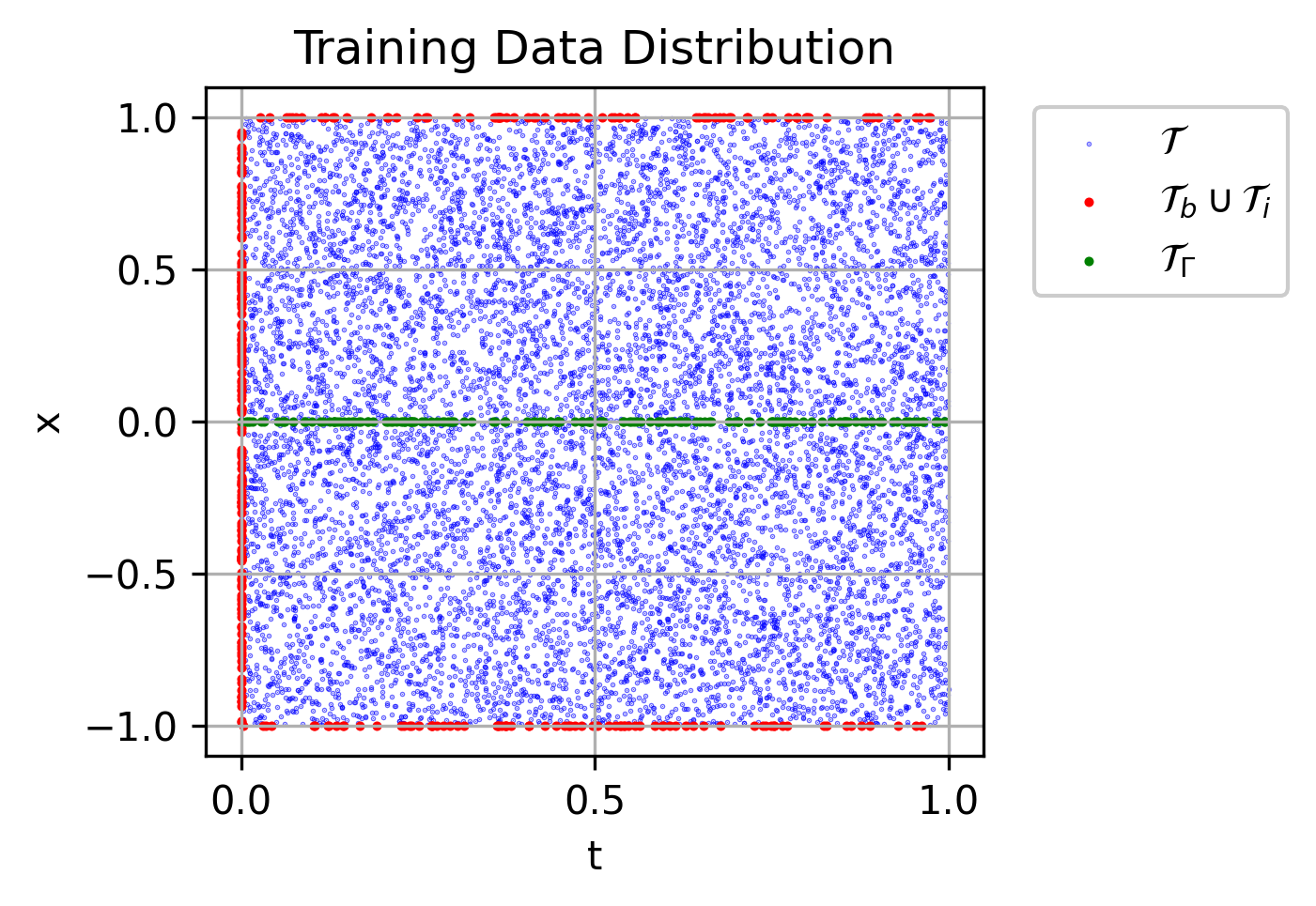}
    \caption{Smooth initial data case: sl-PINN training data distribution.}
    \label{train_data_distribution smooth}
\end{figure}

For both classical PINNs and sl-PINNs, 
we set the neural networks with size $4\times 20$ (3 hidden layers with 20 neurons each) and randomly choose $N=5,000$ training points inside the space-time domain $\Omega_{x} \times \Omega_{t}$, $N_{b}=80$ on the boundary $\partial \Omega_{x}  \times \Omega_{t}$, $N_{0}=80$ on the initial boundary $\Omega_{x} \times \{ t=0 \}$, and $N_{i}=80$ on the interface $\{ x=0 \} \times \Omega_{t}$. The training data distribution is illustrated in Figure \ref{train_data_distribution smooth}. 
For the optimization, we use the Adam optimizer with learning rate $10^{-3}$ for $20,000$ iterations, and then the L-BFGS optimizer with learning rate $1$ for $10,000$ iterations for further training. 
The main reason for using the combination of Adam and L-BFGS is to employ a fast converging (2nd order) optimizer (L-BFGS) after obtaining good initial parameters with the 1st order optimizer (Adam).

\begin{figure}
    \centering
    \begin{subfigure}[b]{0.45\textwidth}
        \centering
        \includegraphics[width=6cm]{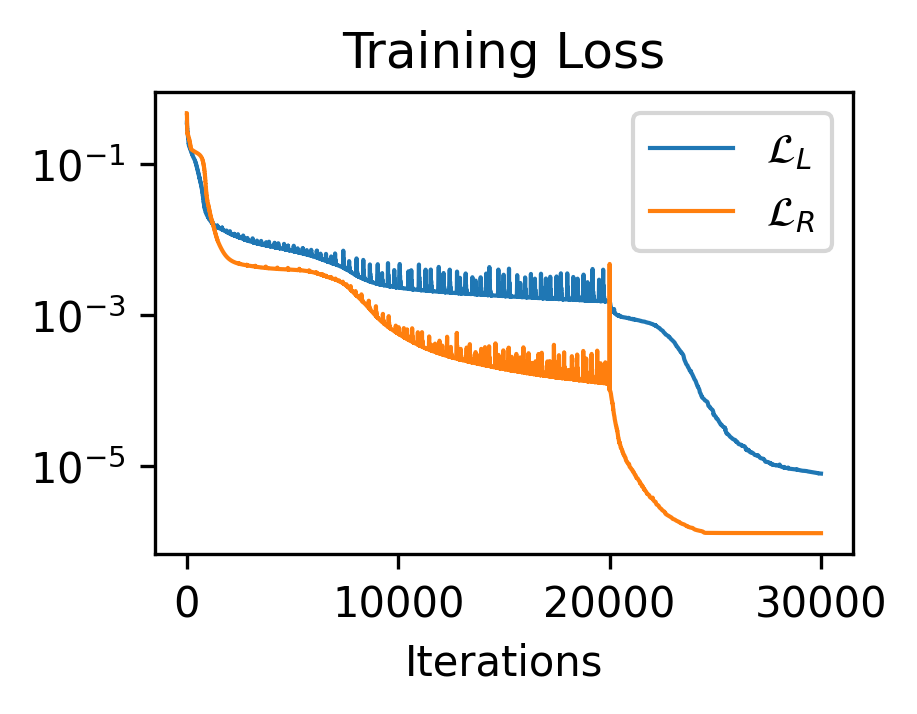}
        \caption{$\ep=10^{-2}/\pi$}
    \end{subfigure}
    \begin{subfigure}[b]{0.45\textwidth}
        \centering
        \includegraphics[width=6cm]{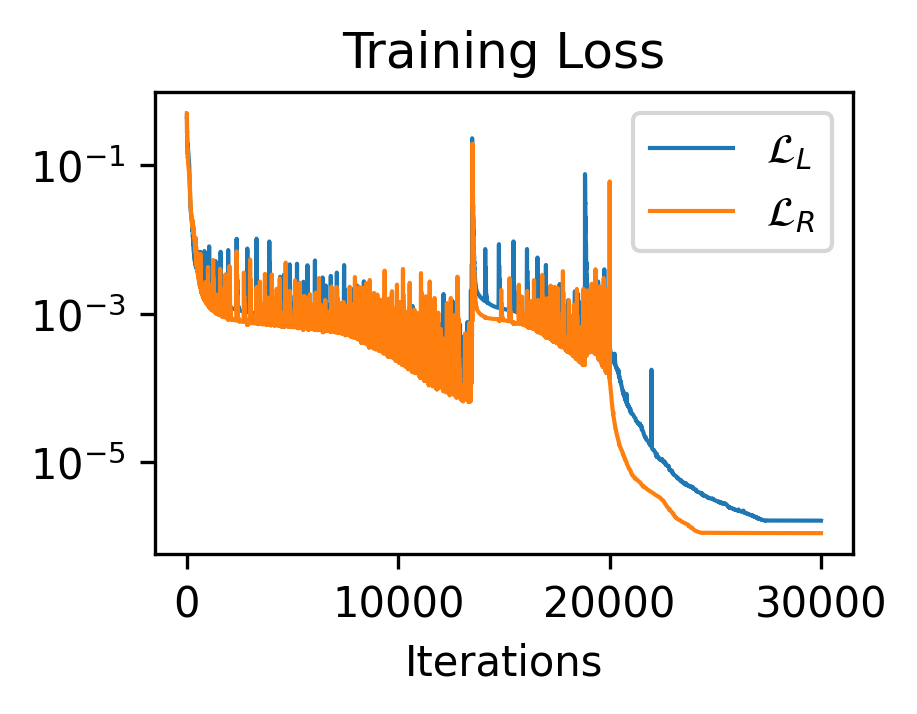}
        \caption{$\ep=10^{-3}/\pi$}
    \end{subfigure}
    \caption{Smooth initial data: sl-PINN's training loss during the training process.}
    \label{smooth train loss slPINN}
\end{figure}

\begin{figure}
    \centering
    \begin{subfigure}[b]{0.45\textwidth}
        \centering
        \includegraphics[width=6cm]{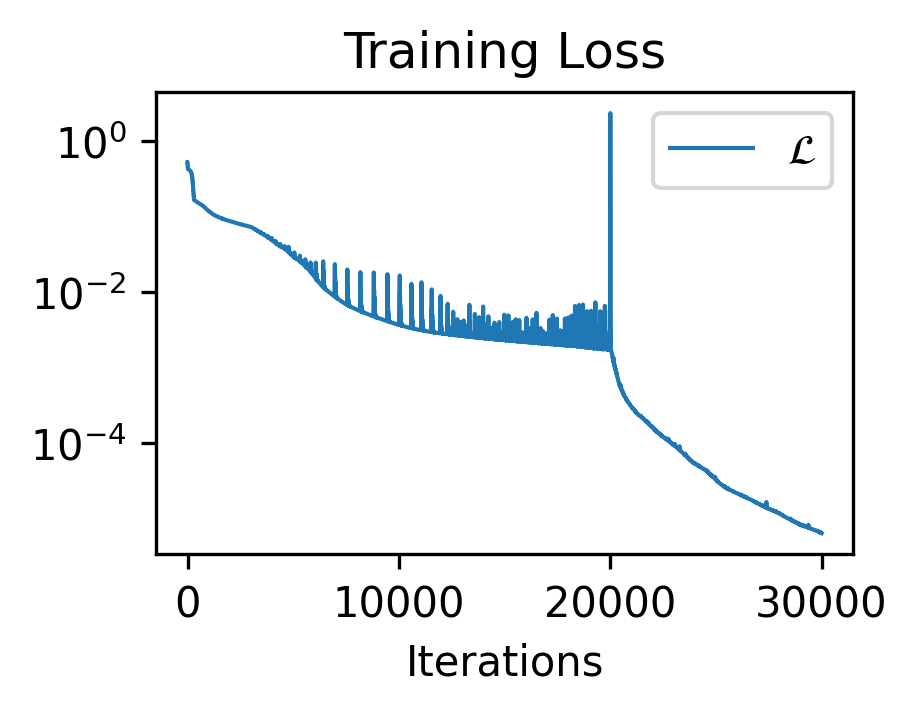}
        \caption{$\ep=10^{-2}/\pi$}
    \end{subfigure}
    \begin{subfigure}[b]{0.45\textwidth}
        \centering
        \includegraphics[width=6cm]{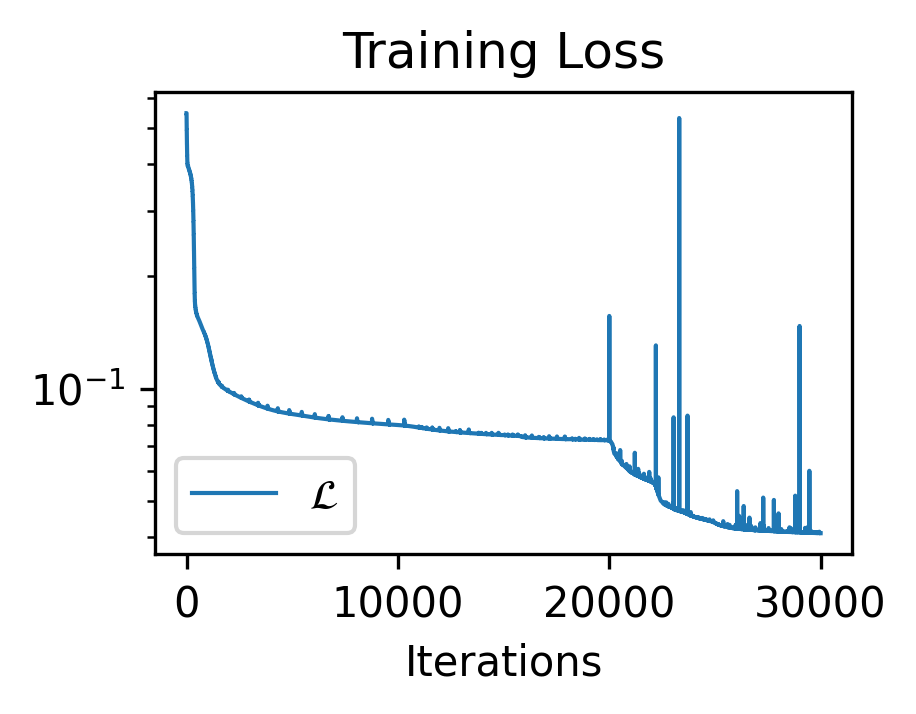}
        \caption{$\ep=10^{-3}/\pi$}
    \end{subfigure}
    \caption{Smooth initial data: PINN's training loss during the training process.}
    \label{smooth train loss PINN}
\end{figure}

Experiments are conducted for both PINNs and sl-PINNs where 
$l(\hat{u})$ in (\ref{e:loss_PINN}) and 
$l(\tilde{u}_*)$, $*=L,R$, in (\ref{loss slPINN}) are defined by 
\begin{equation}\label{e:l_PINN_smooth} 
        l(\hat{u}) 
            = \big| \hat{u}(1, t) - \hat{u}(-1, t) \big|^2
\end{equation}
and
\begin{equation}\label{e:l_slPINN_smooth} 
        l(\tilde{u}_*) 
            = 
                \big| \tilde{u}_R(1, t) - \tilde{u}_L(-1, t) \big|^2, 
        \quad
        \text{for both } 
        * = L, R.
\end{equation}
We repeatedly train the PINNs and sl-PINNs for the Burgers' equations using the viscosity parameter values of $\epsilon=10^{-1}/\pi, \ 10^{-2}/\pi$, and $10^{-3}/\pi$. We use the same computational settings mentioned earlier. Below, we present the training processes for the cases of $\epsilon=10^{-2}/\pi$ and $\epsilon=10^{-3}/\pi$ in Figures \ref{smooth train loss slPINN} and \ref{smooth train loss PINN}.

To evaluate the performance of sl-PINNs, we present the $L^{2}$-error (\ref{L2L2}) and the $L^{\infty}$-error (\ref{Linf}) at a specific time in Tables \ref{table smooth small} and \ref{table smooth large}. We observe that sl-PINN provides accurate predictions for all cases when $\ep=10^{-1}/\pi$, $10^{-2}/\pi$, and $10^{-3}/\pi$ in $L^{2}$-norm, while PINN struggles to perform well when $\ep=10^{-3}/\pi$. Additionally, the $L^{\infty}$-error at a specific time shows that sl-PINN has better accuracy than PINN when the shock occurs as time progresses.

\begin{figure}
    \centering
    \begin{subfigure}[b]{0.4\textwidth}
        \centering
        \includegraphics[width=\textwidth]{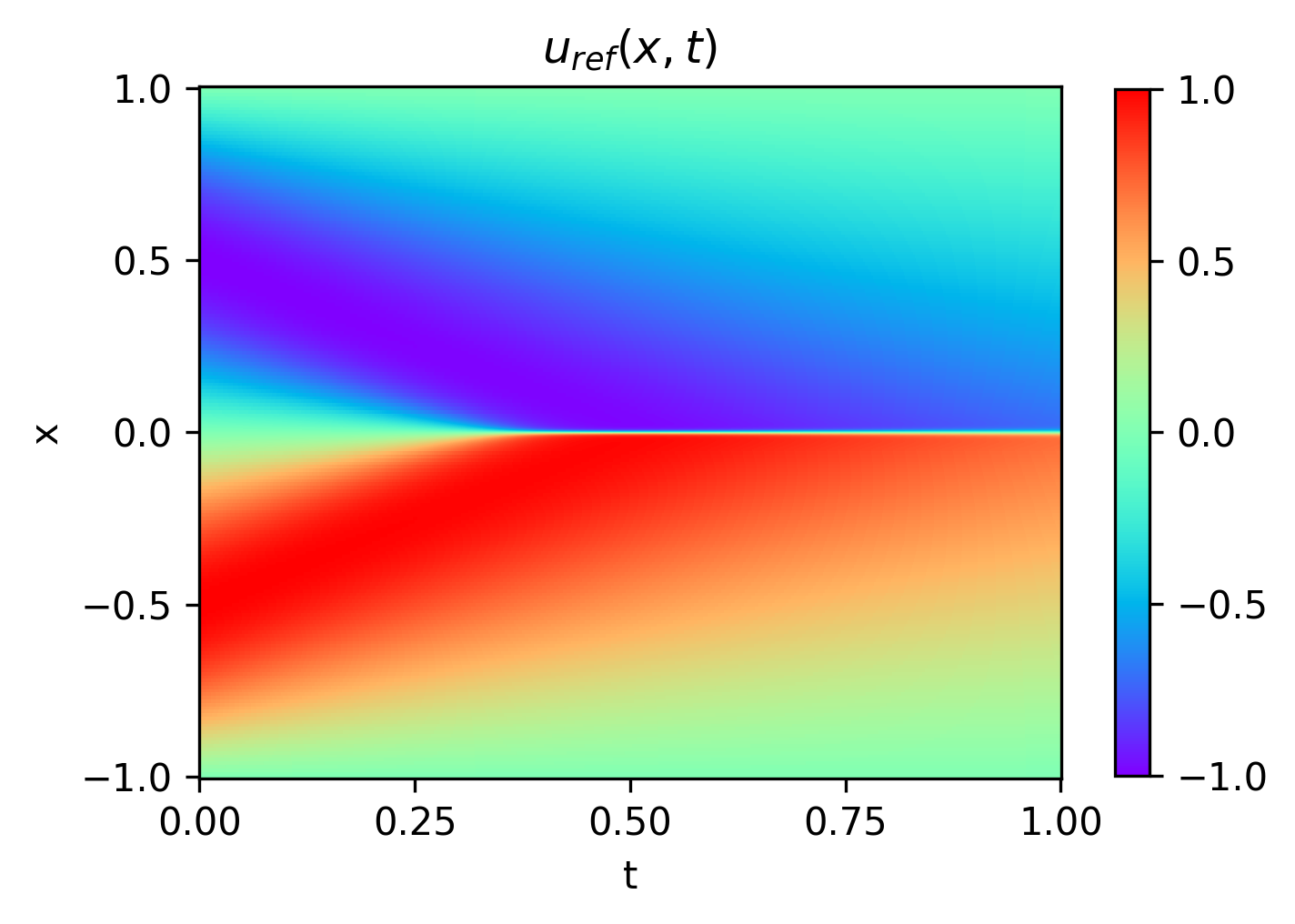}
        \caption{reference solution: $\ep=10^{-2}/\pi$}
    \end{subfigure}
    
    \bigskip
    
    \begin{subfigure}[b]{0.8\textwidth}
        \centering
        \includegraphics[width=\textwidth]{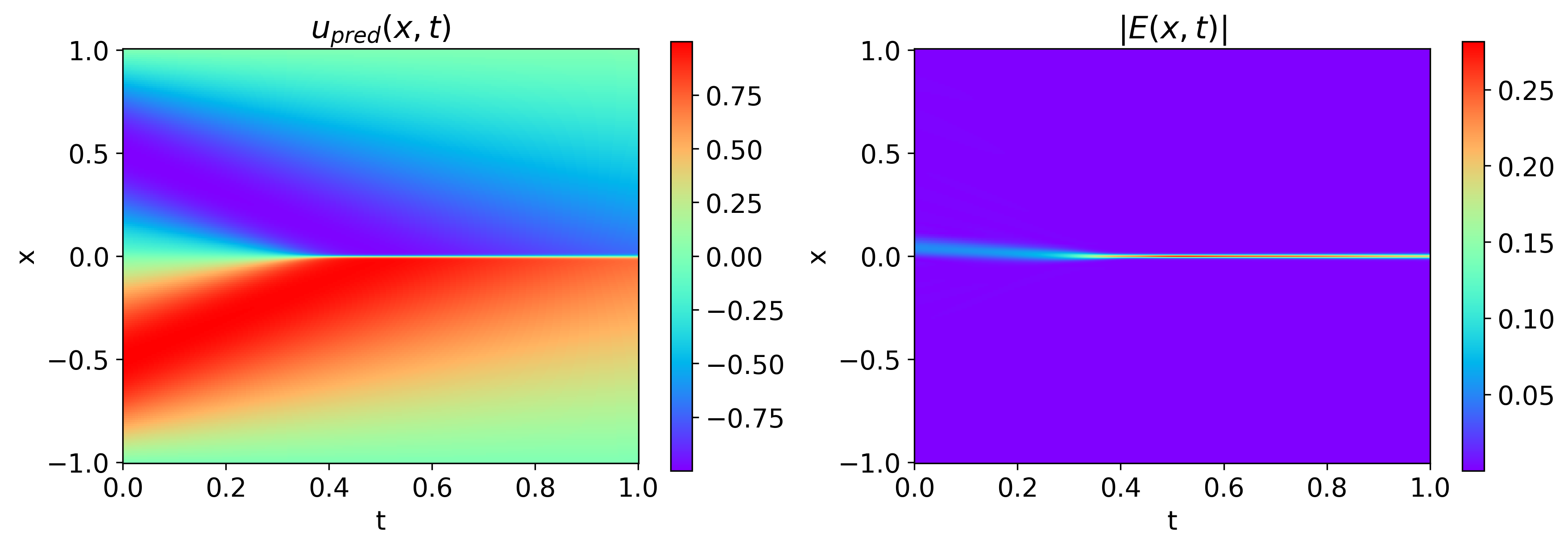}
        \caption{PINN prediction (left) and absolute pointwise error (right)}
    \end{subfigure}
    
    \bigskip
    
    \begin{subfigure}[b]{0.8\textwidth}
        \centering
        \includegraphics[width=\textwidth]{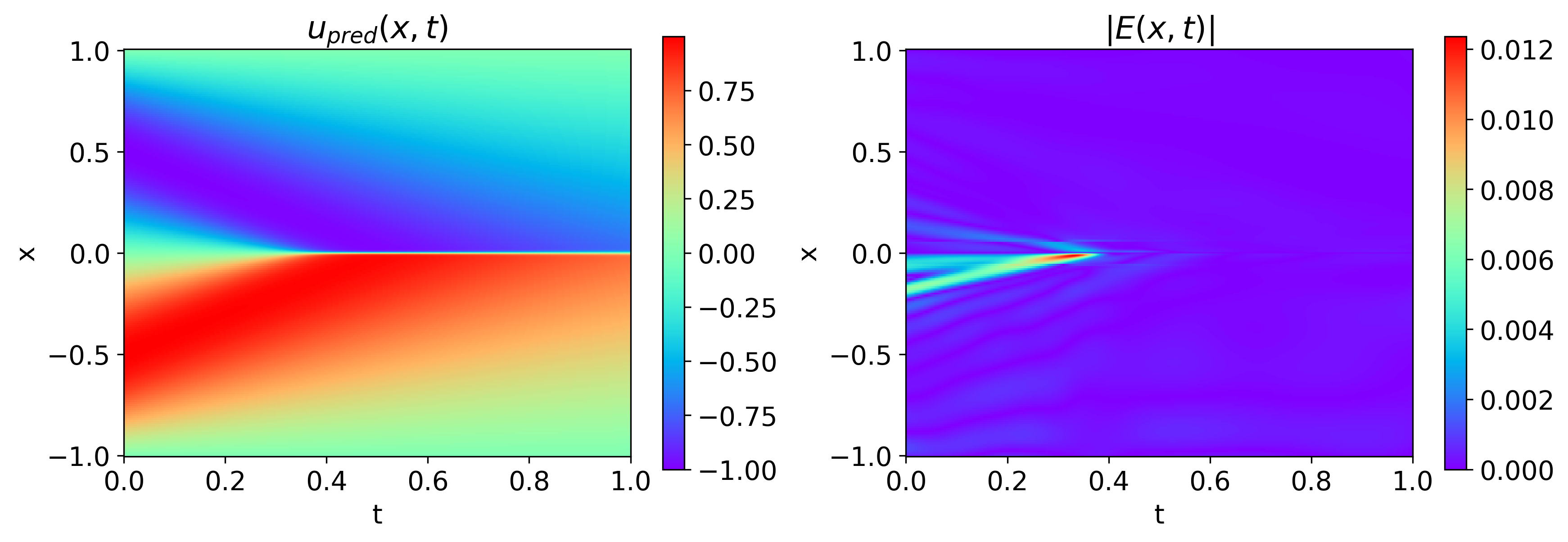}
        \caption{sl-PINN prediction (left) and absolute pointwise error (right)}
    \end{subfigure}
    \caption{Smooth initial data case when $\ep=10^{-2}/\pi$: PINN vs sl-PINN solution plots with small NN size $4\times 20$.}
    \label{plot2d smooth 1}
\end{figure}

\begin{figure}
    \centering
    \begin{subfigure}[b]{0.4\textwidth}
        \centering
        \includegraphics[width=\textwidth]{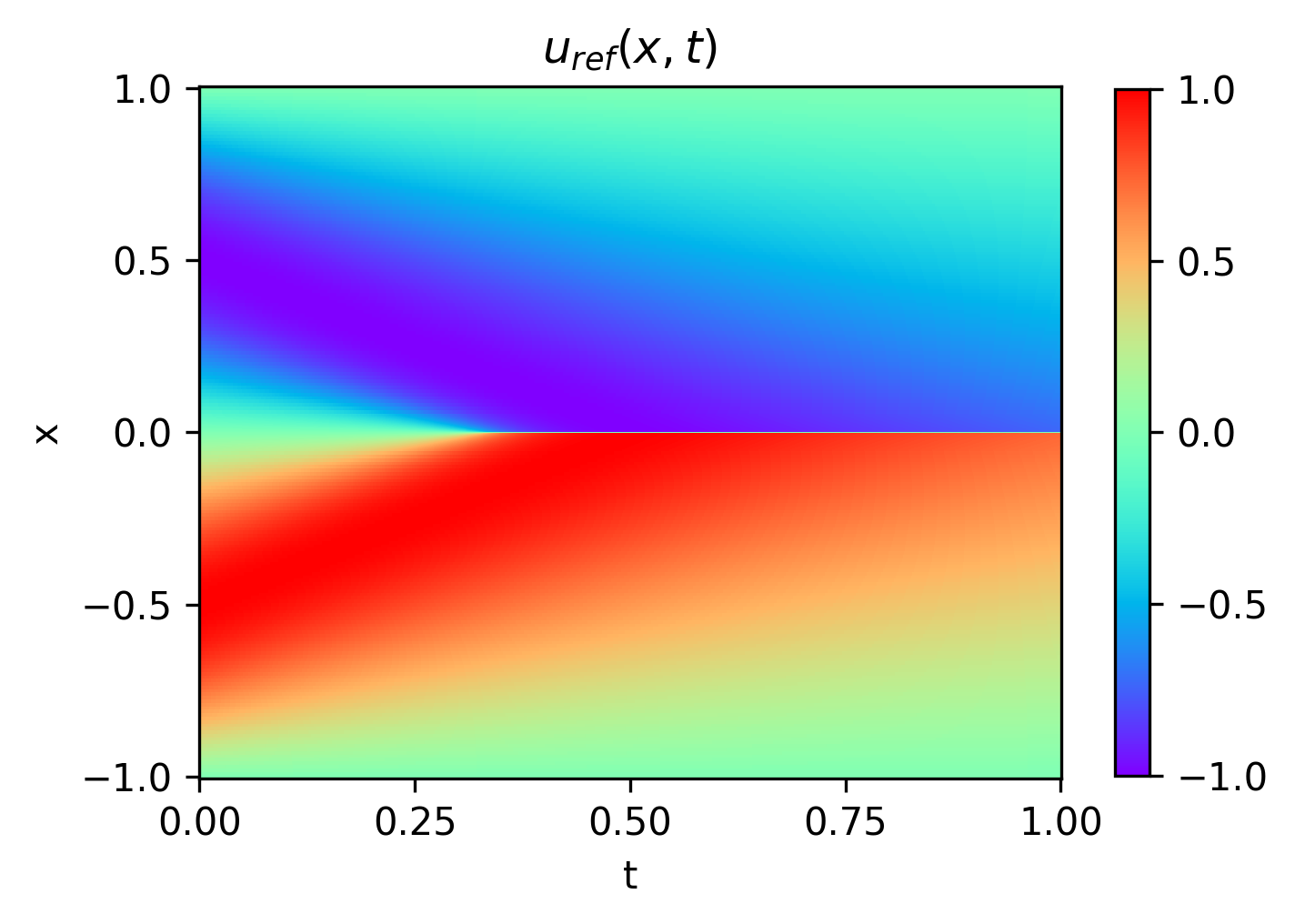}
        \caption{reference solution: $\ep=10^{-3}/\pi$}
    \end{subfigure}
    
    \bigskip
    
    \begin{subfigure}[b]{0.8\textwidth}
        \centering
        \includegraphics[width=\textwidth]{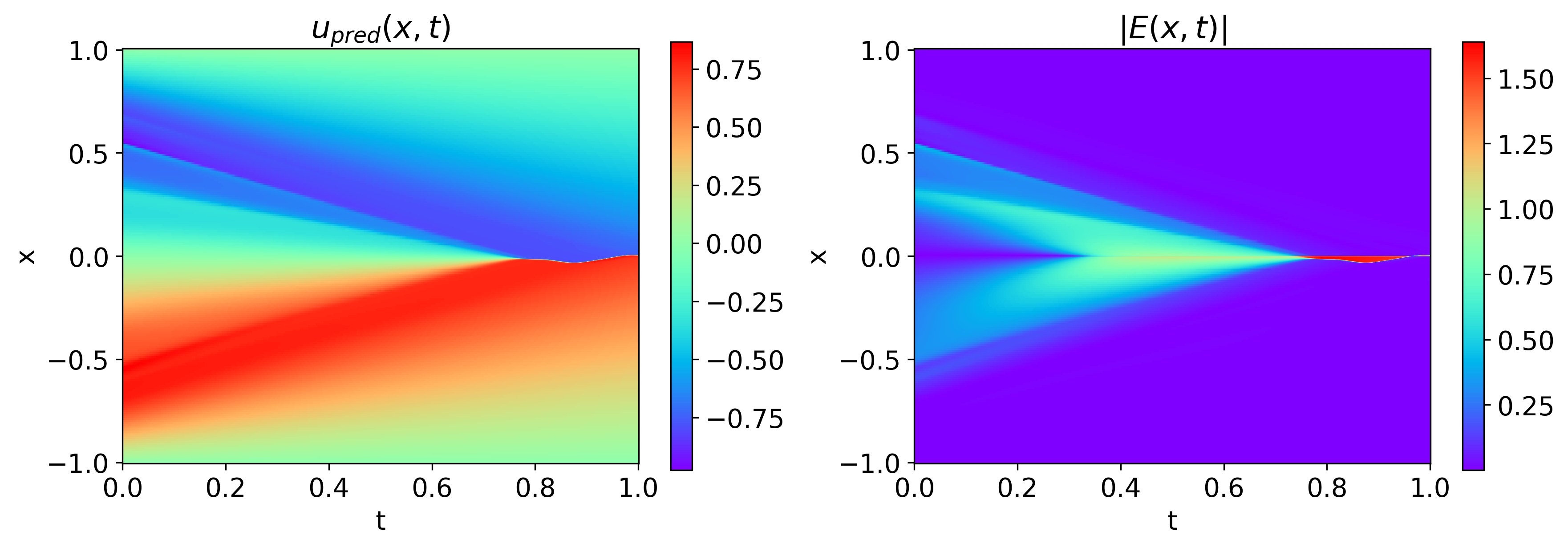}
        \caption{PINN prediction (left) and absolute pointwise error (right)}
    \end{subfigure}
    
    \bigskip
    
    \begin{subfigure}[b]{0.8\textwidth}
        \centering
        \includegraphics[width=\textwidth]{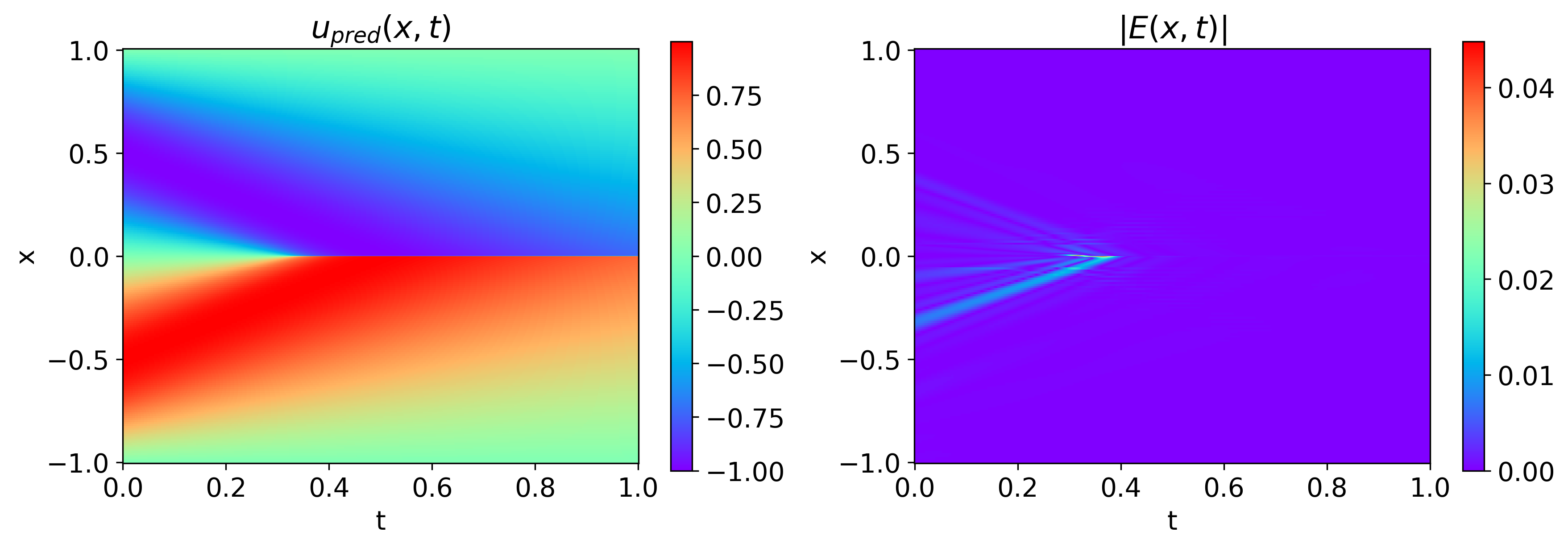}
        \caption{sl-PINN prediction (left) and absolute pointwise error (right)}
    \end{subfigure}
    \caption{Smooth initial data case when $\ep=10^{-3}/\pi$: PINN vs sl-PINN solution plots with small NN size $4\times 20$.}
    \label{plot2d smooth 2}
\end{figure}

\begin{figure}
    \centering
    \begin{subfigure}[b]{\textwidth}
        \centering
        \includegraphics[width=\textwidth]{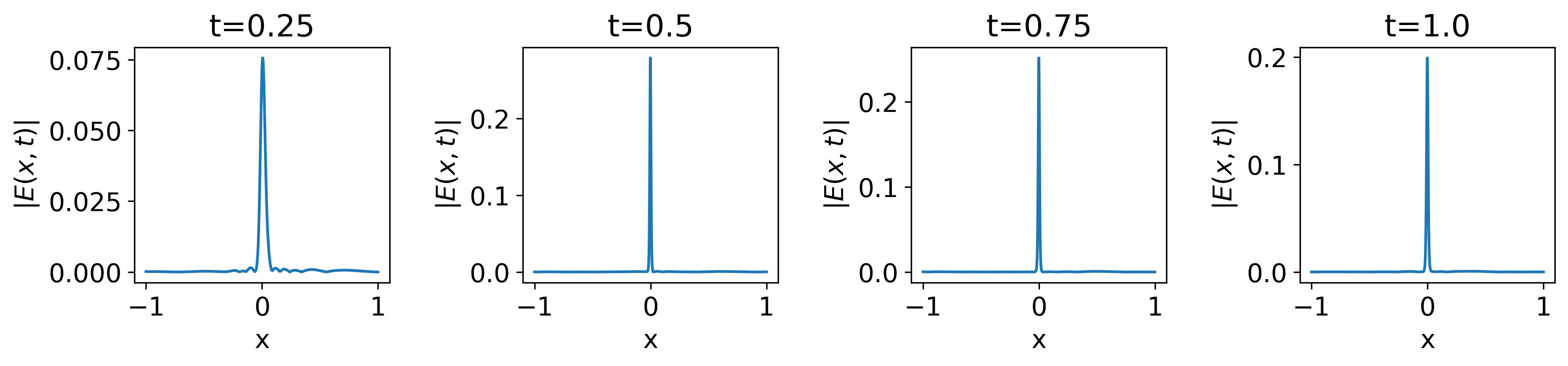}
        \caption{Errors of PINN prediction at specific time.}
    \end{subfigure}
    
    \bigskip
    
    \begin{subfigure}[b]{\textwidth}
        \centering
        \includegraphics[width=\textwidth]{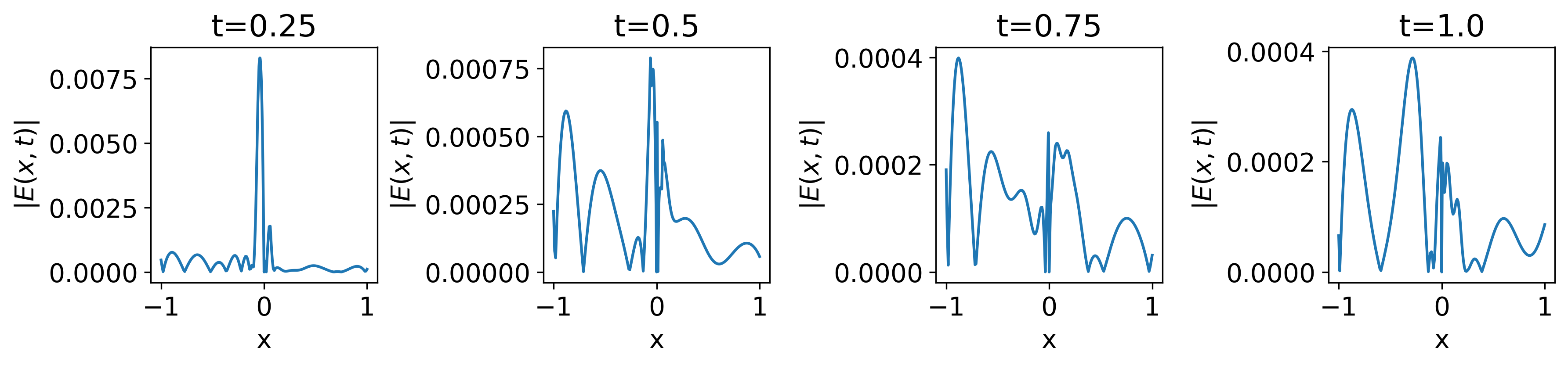}
        \caption{Errors of sl-PINN prediction at specific time.}
    \end{subfigure}
    \caption{Smooth initial data case when $\ep=10^{-2}/\pi$: PINN vs sl-PINN solution error at specific time (small NN size $4\times 20$).}
    \label{plot1d smooth 1}
\end{figure}

\begin{figure}
    \centering
    \begin{subfigure}[b]{\textwidth}
        \centering
        \includegraphics[width=\textwidth]{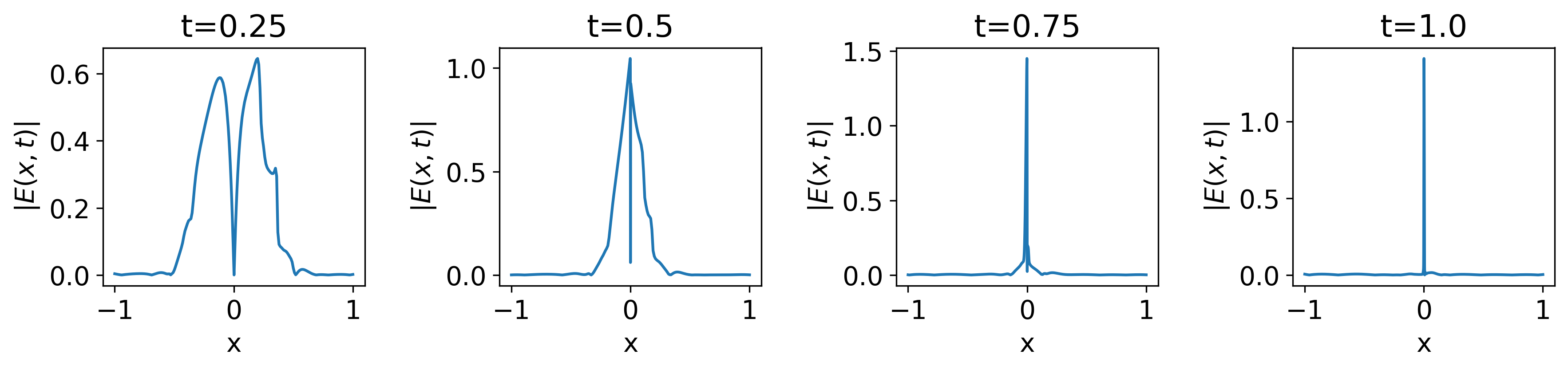}
        \caption{Errors of PINN prediction at specific time.}
    \end{subfigure}
    
    \bigskip
    
    \begin{subfigure}[b]{\textwidth}
        \centering
        \includegraphics[width=\textwidth]{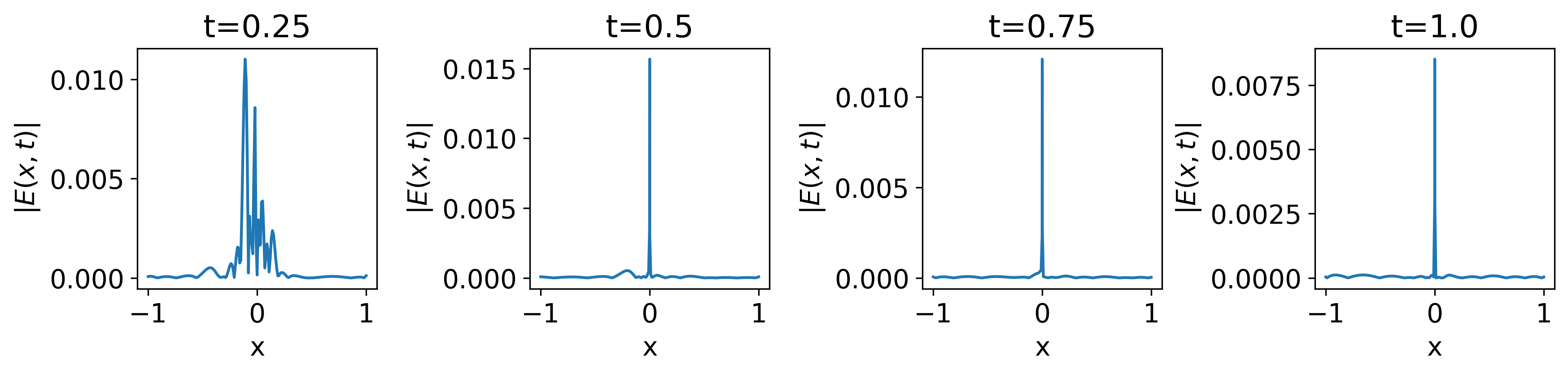}
        \caption{Errors of sl-PINN prediction at specific time.}
    \end{subfigure}
    \caption{Smooth initial data case when $\ep=10^{-3}/\pi$: PINN vs sl-PINN solution error at specific time (small NN size $4\times 20$).}
    \label{plot1d smooth 2}
\end{figure}

The figures below show the profiles of the PINN and sl-PINN predictions. In particular, Figure \ref{plot2d smooth 1} illustrates the case when $\ep=10^{-2}/\pi$, and Figure \ref{plot2d smooth 2} illustrates the case when $\ep=10^{-3}/\pi$. It is observed that for $\ep=10^{-2}/\pi$, the sl-PINN method is successful in resolving significant errors near the shock that occurred in the PINN prediction. Additionally, when $\ep=10^{-3}/\pi$, the sl-PINN method successfully predicts the solution with a stiffer interior layer which the PINN fails to predict. Furthermore, Figure \ref{plot1d smooth 1} shows the temporal snapshots of the errors of PINN and sl-PINN for $\ep=10^{-2}/\pi$, while Figure \ref{plot1d smooth 2} shows the same for $\ep=10^{-3}/\pi$.

In a study addressing the use of Physics Informed Neural Networks (PINNs) to solve the Burgers equation (referenced as \cite{RPK19}), it was found that increasing the number of layers in the neural network leads to improved predictive accuracy.  
To further test this, a deeper neural network with a size of $9 \times 20$ (8 hidden layers with 20 neurons each) was employed, along with a doubling of the training points to $N=10,000$, $N_{b}=160$, and $N_{0}=160$. As indicated in Table \ref{table smooth large}, when $\ep=10^{-2}/\pi$, the accuracy of PINN increased. However, there was no improvement in accuracy for the stiffer case when $\ep=10^{-3}/\pi$. Interestingly, despite using a larger neural network, the singular-layer PINNs (sl-PINNs) performe similarly to using smaller neural networks. Consequently, it was concluded that utilizing small neural networks of size $4\times 20$ for sl-PINNs is sufficient to achieve a good approximation, thereby increasing computational efficiency.  

Additionally, it's important to note that a two-scale neural networks learning method for the Burgers equation with small $\ep>0$ was proposed in a recent study (referenced as \cite{ZYZK24}). This approach improved the accuracy of the predicted solution by incorporating stretch variables as additional inputs to the neural network. However, this method  requires appropriate initialization from pre-trained models. Sequential training is necessary to obtain the pre-trained models for $\ep=10^{-1}/\pi$ and $\ep=10^{-2}/\pi$ in order to achieve better initialization of the neural network for training the most extreme case, $\ep=10^{-3}/\pi$.  
In contrast, our sl-PINN method only employs two neural networks for the case $\ep=10^{-3}/\pi$ and it still achieves better accuracy.


\begin{table}[]
\centering
\begin{subtable}{\textwidth}
\centering
\begin{tabular}{|l|r|r|}
\hline
$\| E \|_{L^{2}_{t}( L^{2}_{x}(\Omega))}$ & \multicolumn{1}{l|}{PINNs} & \multicolumn{1}{l|}{sl-PINNs} \\ \hline
$\ep=10^{-1}/\pi$                         & 4.30E-04                   & 4.28E-04                      \\
$\ep=10^{-2}/\pi$                         & 6.02E-02                   & 2.13E-03                      \\
$\ep=10^{-3}/\pi$                         & 5.50E-01                   & 5.57E-03                      \\ \hline
\end{tabular}
\bigskip
\end{subtable}
\begin{subtable}{\textwidth}
\centering
\begin{tabular}{|l|rr|rr|}
\hline
                                   & \multicolumn{2}{c|}{PINNs}                                                      & \multicolumn{2}{c|}{sl-PINNs}                                                   \\ \hline
$\| E \|_{L^{\infty}_{x}(\Omega)}$ & \multicolumn{1}{l|}{$\ep=10^{-2}/\pi$} & \multicolumn{1}{l|}{$\ep=10^{-3}/\pi$} & \multicolumn{1}{l|}{$\ep=10^{-2}/\pi$} & \multicolumn{1}{l|}{$\ep=10^{-3}/\pi$} \\ \hline
$t = 0.25$                         & \multicolumn{1}{r|}{7.56E-02}          & 6.45E-01                               & \multicolumn{1}{r|}{8.31E-03}          & 1.10E-02                               \\
$t = 0.5$                          & \multicolumn{1}{r|}{2.79E-01}          & 1.05E+00                               & \multicolumn{1}{r|}{7.89E-04}          & 1.57E-02                               \\
$t = 0.75$                         & \multicolumn{1}{r|}{2.52E-01}          & 1.45E+00                               & \multicolumn{1}{r|}{3.99E-04}          & 1.21E-02                               \\
$t = 1.0$                          & \multicolumn{1}{r|}{1.99E-01}          & 1.41E+00                               & \multicolumn{1}{r|}{3.88E-04}          & 8.51E-03                               \\ \hline
\end{tabular}
\end{subtable}
\caption{Error comparison  between PINN and sl-PINN: $L^{2}$-errors (top) and $L^{\infty}$-errors (bottom) with small NN size: $4\times 20$.}
\label{table smooth small}
\end{table}

\begin{table}[]
\centering
\begin{subtable}{\textwidth}
\centering
\begin{tabular}{|l|r|r|}
\hline
$\| E \|_{L^{2}_{t}( L^{2}_{x}(\Omega))}$ & \multicolumn{1}{l|}{PINNs} & \multicolumn{1}{l|}{sl-PINNs} \\ \hline
$\ep=10^{-1}/\pi$                         & 4.23E-04                   & 2.15E-04                      \\
$\ep=10^{-2}/\pi$                         & 9.43E-04                   & 7.10E-04                      \\
$\ep=10^{-3}/\pi$                         & 4.83E-01                   & 4.25E-03                      \\ \hline
\end{tabular}
\bigskip
\end{subtable}
\begin{subtable}{\textwidth}
\centering
\begin{tabular}{|l|rr|rr|}
\hline
                                   & \multicolumn{2}{c|}{PINNs}                                                      & \multicolumn{2}{c|}{sl-PINNs}                                                   \\ \hline
$\| E \|_{L^{\infty}_{x}(\Omega)}$ & \multicolumn{1}{l|}{$\ep=10^{-2}/\pi$} & \multicolumn{1}{l|}{$\ep=10^{-3}/\pi$} & \multicolumn{1}{l|}{$\ep=10^{-2}/\pi$} & \multicolumn{1}{l|}{$\ep=10^{-3}/\pi$} \\ \hline
$t = 0.25$                         & \multicolumn{1}{r|}{1.51E-03}          & 5.91E-01                               & \multicolumn{1}{r|}{1.38E-03}          & 1.02E-02                               \\
$t = 0.5$                          & \multicolumn{1}{r|}{2.70E-03}          & 1.04E+00                               & \multicolumn{1}{r|}{3.62E-04}          & 1.57E-02                               \\
$t = 0.75$                         & \multicolumn{1}{r|}{3.26E-03}          & 1.52E+00                               & \multicolumn{1}{r|}{3.33E-04}          & 1.21E-02                               \\
$t = 1.0$                          & \multicolumn{1}{r|}{2.03E-03}          & 7.65E-01                               & \multicolumn{1}{r|}{2.14E-04}          & 8.51E-03                               \\ \hline
\end{tabular}
\end{subtable}
\caption{Error comparison  between PINN and sl-PINN: $L^{2}$-errors (top) and $L^{\infty}$-errors (bottom) with large NN size: $9\times 20$.}
\label{table smooth large}
\end{table}

\subsection{Non-smooth initial data} 

We consider the Burgers' equation (\ref{Burgers}) 
supplemented with 
the boundary condition $(ii)$ in the form,     
\begin{equation}\label{non smooth BC -gen}
    \lim_{x \rightarrow - \infty} u^\ep  = u_L, \quad
    \lim_{x \rightarrow \infty} u^\ep = -u_R, \quad
    u_L, \, u_R > 0, 
\end{equation}
and  the initial condition, 
\begin{equation}\label{non smooth IC}
 u_{0}(x) =
    \begin{cases}
       u_L, &  x < 0, \\
       -u_R, &  x > 0.
    \end{cases}
\end{equation}

Using the Cole-Hopf transformation; see, e.g, \cite{Cole, Hopf, Peter2014}, 
we find the explicit expression of solution to the viscous Burgers' equation  (\ref{Burgers}) with (\ref{non smooth BC -gen}) and (\ref{non smooth IC}) as 
\begin{align}\label{e:exact_nonsmooth}
    u^\ep (x,t) 
        = 
            u_L 
            - 
            \frac{u_L + u_R}
                {1 
                + 
                \exp\bigg(- \dfrac{u_L + u_R}{2\ep}(x - ct)\bigg) 
                \text{erfc} \left( \dfrac{x- u_L t}{2\sqrt{\ep t}} \right) \big
                / 
                \text{erfc} \left( \dfrac{- x- u_R t}{2\sqrt{\ep t}} \right) }.
\end{align}
Here $c=(u_L - u_R)/2$ and $\text{erfc}(\cdot)$ is the complementary error function on $(-\infty, \infty )$
defined by 
\begin{equation}
    \text{erfc}(z) 
        = 1 - \text{erf}(z)
        = 1 - \dfrac{2}{\sqrt{\pi}} \int_0^z e^{-t^2} \, dt.
\end{equation}
Note that accurately evaluating the values of the complementary error function is a significantly challenging task. 
Therefore, to obtain the reference solution to (\ref{Burgers}) with (\ref{non smooth BC -gen}) and (\ref{non smooth IC}) 
in the upcoming numerical computations below, 
we use the high-resolution finite difference method as described.



The inviscid limit problem (\ref{Burgers_inviscid}) with the boundary and initial data (\ref{non smooth BC -gen}) and (\ref{non smooth IC}), known as the Riemann problem, 
has a unique weak solution (see, e.g., \cite{Howison, Peter2014}) given by
\begin{equation*}
    u^{0}(x,t) = 
    \begin{cases}
        u_L, &  x < st, \\
        - u_R, &  x > st,
    \end{cases}
\end{equation*}
where $s$ is the speed of the shock wave, determined by the Rankine-Hugoniot jump condition, written as 
\begin{align*}
    s(t) 
    = \alpha^\prime(t)
    = - \dfrac{1}{2} \dfrac{ u_R^2 - u_L^2}{ u_R + u_L} = \dfrac{1}{2} (u_L - u_R).
\end{align*}
Integrating the equation above, we find the shock curve 
\begin{align*}
    \alpha(t) = \frac{1}{2}(u_L - u_R) t, \quad t \geq0.
\end{align*}

In the following experiments, we consider two examples: the first is when $u_L = u_R = 1$ for a steady shock located at $x=\alpha(t) \equiv 0$, and the second is when $u_L=1$ and $u_R=1/2$ for a moving shock located at $x=\alpha(t)=t/4 > 0$. 
Experiments for both examples are conducted for both PINNs and sl-PINNs especially when the viscosity is small 
and hence the value of $u^\ep$ at $x=-1$ (or at $x= 1$) is exponentially close to the boundary value $u_L$ (or $-u_R$) at least for the time of computations. 
Thanks to this property, for the sake of computational convenience, we set  
$l(\hat{u})$ in (\ref{e:loss_PINN}) and 
$l(\tilde{u}_*)$, $*=L,R$, in (\ref{loss slPINN}) as 
\begin{equation}\label{e:l_PINN_nonsmooth} 
        l(\hat{u}) 
            = 
                \big| \hat{u}(-1, t) - u_L \big|^2 
                +
                \big| \hat{u}(1, t) + u_R \big|^2,
\end{equation}
and
\begin{equation}\label{e:l_slPINN_nonsmooth} 
        l(\tilde{u}_L) 
            = \big| \tilde{u}_L(-1, t) - u_L\big|^2, 
            \qquad
        l(\tilde{u}_R) 
            = \big| \tilde{u}_R(1, t) + u_R \big|^2.
\end{equation}
We repeatedly train the PINNs and sl-PINNs for the Burgers' equations using the viscosity parameter values of $\epsilon=10^{-1}/\pi, \ 10^{-2}/\pi$, and $10^{-3}/\pi$.


{\bf \red }

\begin{figure}
    \centering
    \begin{subfigure}[b]{0.45\textwidth}
        \centering
        \includegraphics[width=8cm]{figures/plot_traindata_u0sin.png}
        \caption{steady shock case}
    \end{subfigure}
    \hfill
    \begin{subfigure}[b]{0.45\textwidth}
        \centering
        \includegraphics[width=8cm]{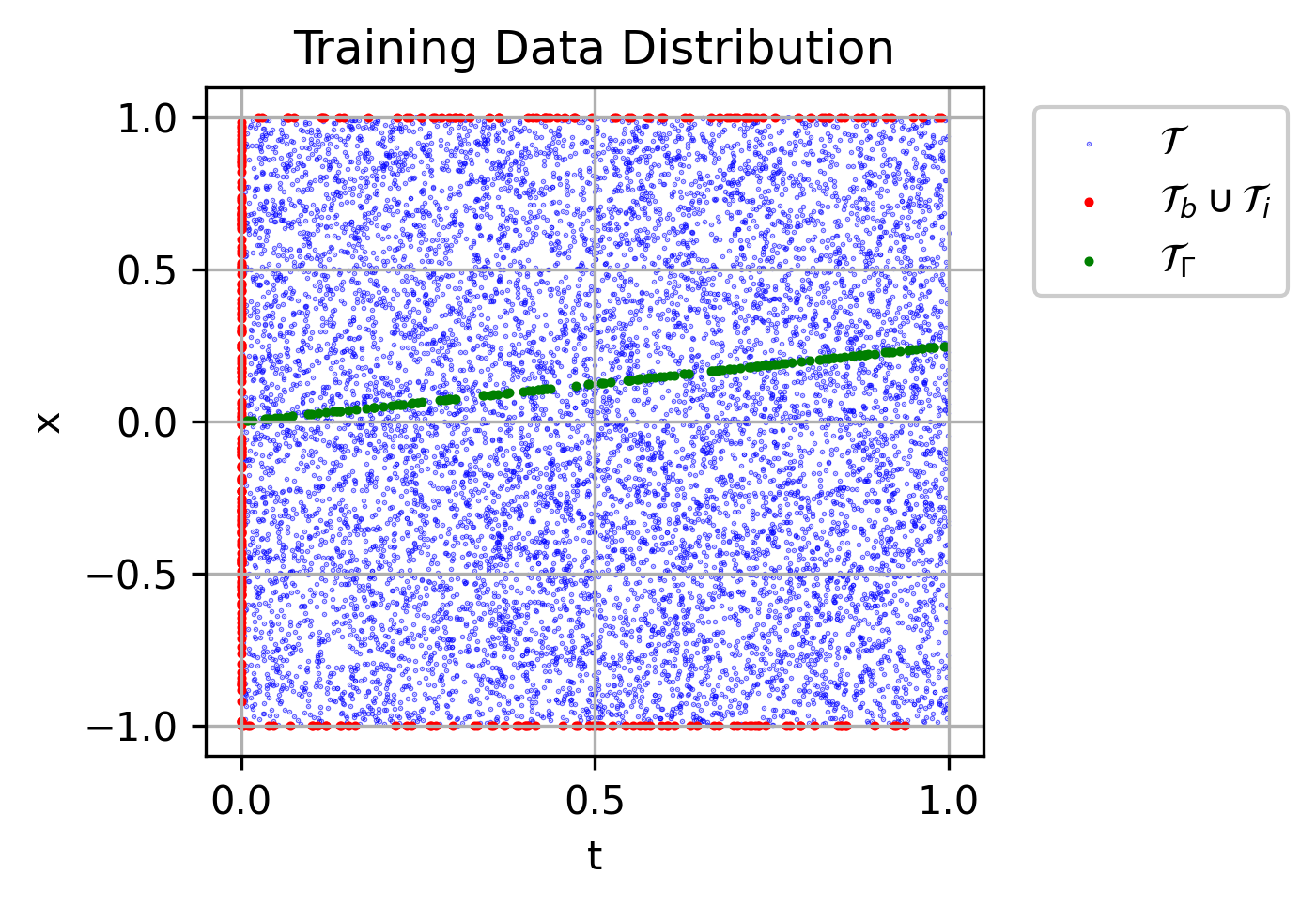}
        \caption{moving shock case}
    \end{subfigure}
    \caption{Nonsmooth initial data: sl-PINN training data distribution.}
    \label{train_data_distribution nonsmooth}
\end{figure}

As for the previous one with smooth data in Section \ref{S.smooth_data}, 
here we use neural networks with a size of $4\times 20$ (3 hidden layers with 20 neurons each). We randomly select $N=5,000$ training points within the space-time domain $\Omega_{x} \times \Omega_{t}$, $N_{b}=80$ on the boundary $\partial \Omega_{x}  \times \Omega_{t}$, $N_{0}=80$ on the initial boundary $\Omega_{x} \times \{ t=0 \}$, and $N_{i}=80$ on the interface $\{ x=\alpha(t) \} \times \Omega_{t}$. The distributions of training data for both steady shock case and moving shock case are shown in Figure \ref{train_data_distribution nonsmooth}. 

\begin{figure}
    \centering
    \begin{subfigure}[b]{0.45\textwidth}
        \centering
        \includegraphics[width=6cm]{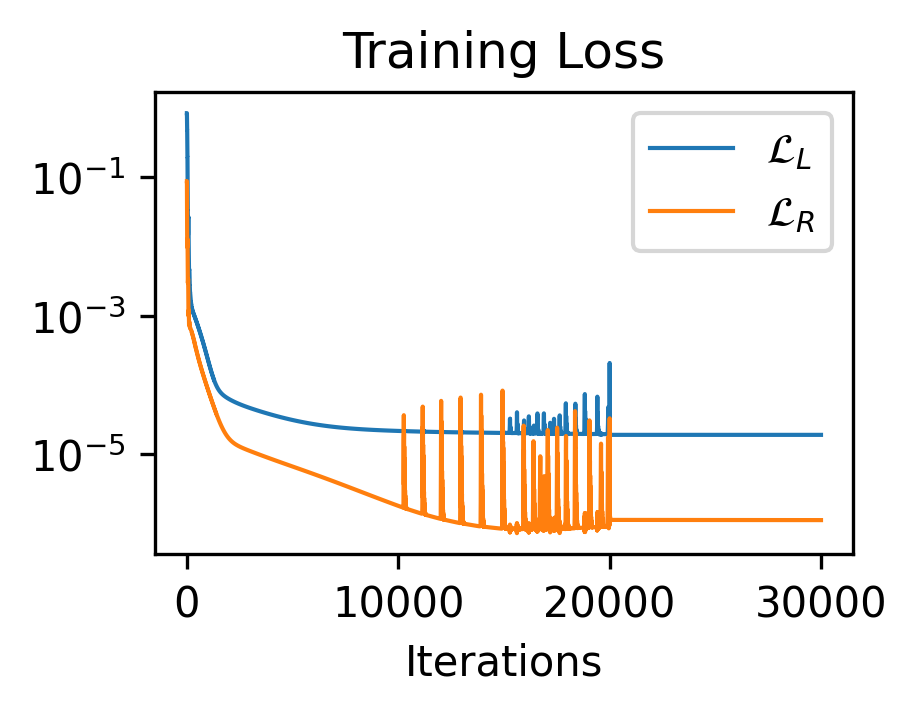}
        \caption{$\ep=1/500$}
    \end{subfigure}
    \begin{subfigure}[b]{0.45\textwidth}
        \centering
        \includegraphics[width=6cm]{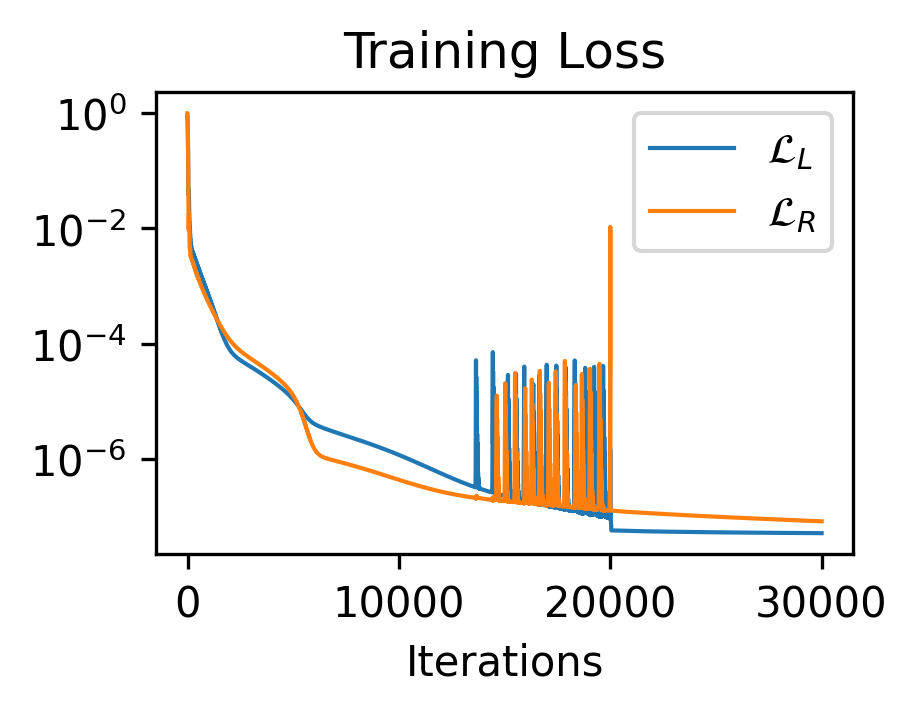}
        \caption{$\ep=1/10000$}
    \end{subfigure}
    \caption{Riemann problem moving shock case: sl-PINN's training loss during the training process.}
    \label{nonsmooth train loss slPINN}
\end{figure}

\begin{figure}
    \centering
    \begin{subfigure}[b]{0.45\textwidth}
        \centering
        \includegraphics[width=6cm]{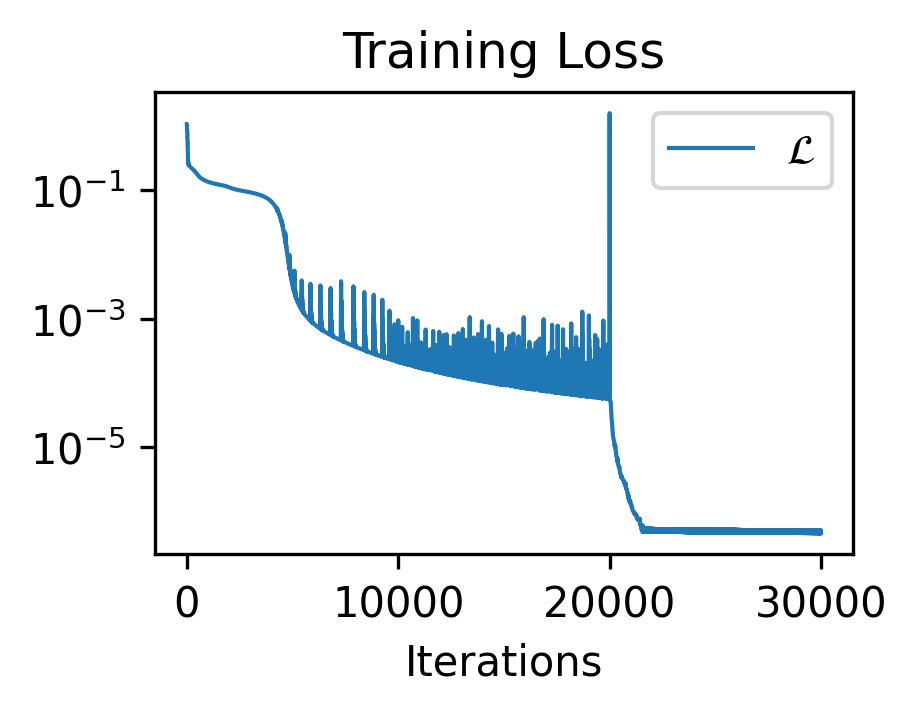}
        \caption{$\ep=1/500$}
    \end{subfigure}
    \begin{subfigure}[b]{0.45\textwidth}
        \centering
        \includegraphics[width=6cm]{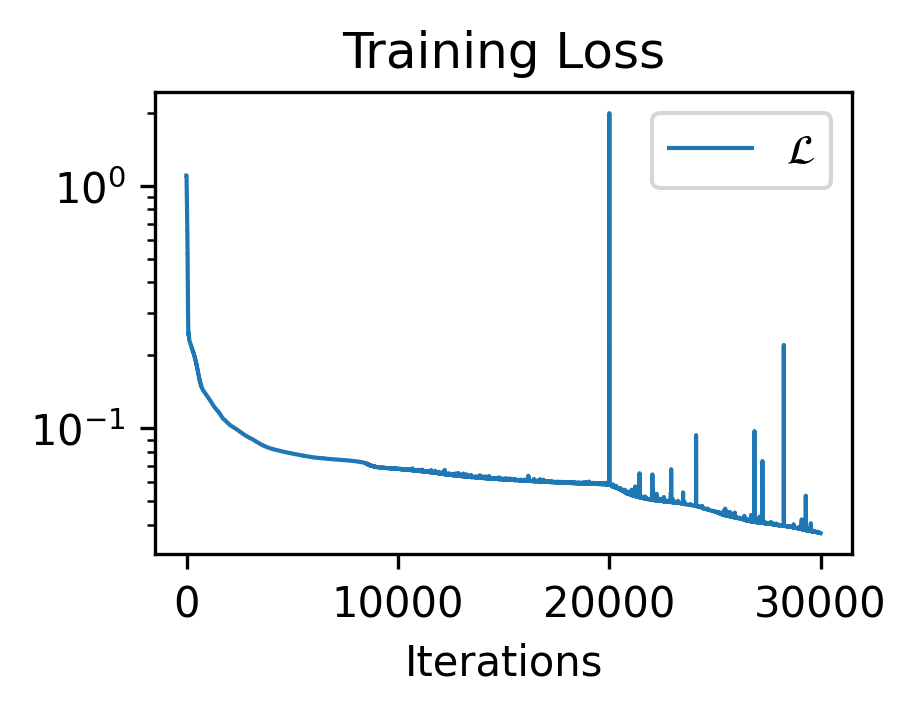}
        \caption{$\ep=1/10000$}
    \end{subfigure}
    \caption{Riemann problem moving shock case: PINN's training loss during the training process.}
    \label{nonsmooth train loss PINN}
\end{figure}

We use the same optimization strategy as for the previous example, which involves adopting the Adam optimizer with a learning rate of $10^{-3}$ for the first 20,000 iteration steps, followed by the L-BFGS optimizer with a learning rate of $1$ for the next 10,000 iterations. We train the cases for $\ep=1/500, \ 1/1000, \ 1/5000,$ and $1/10000$ with the same settings as described. The training losses for the moving shock case for both sl-PINN and PINN when $\ep=1/500$ and $\ep=1/10000$ are presented in Figures \ref{nonsmooth train loss slPINN} and \ref{nonsmooth train loss PINN} respectively.

\begin{figure}
    \centering
    \begin{subfigure}[b]{0.4\textwidth}
        \centering
        \includegraphics[width=\textwidth]{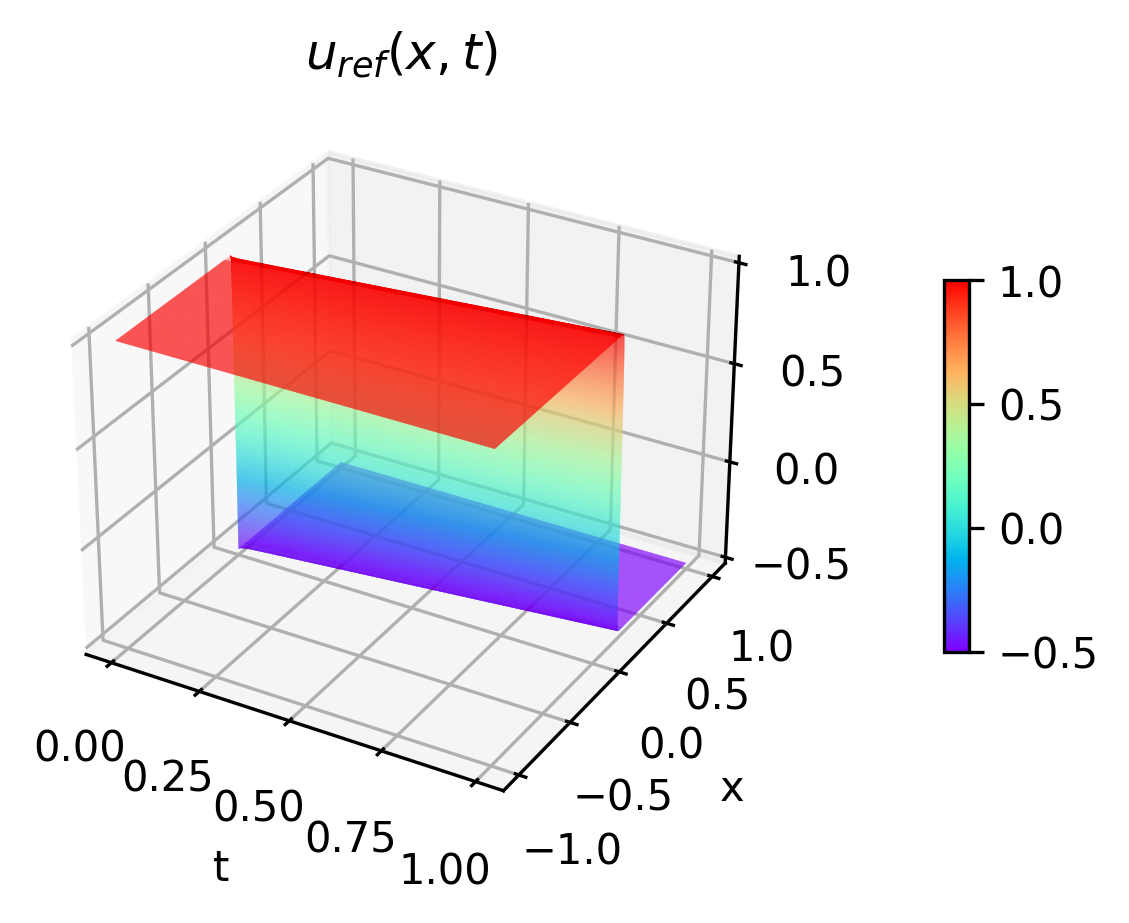}
        \caption{reference solution $\ep=1/500$}
    \end{subfigure}
    
    \bigskip
    
    \begin{subfigure}[b]{0.8\textwidth}
        \centering
        \includegraphics[width=\textwidth]{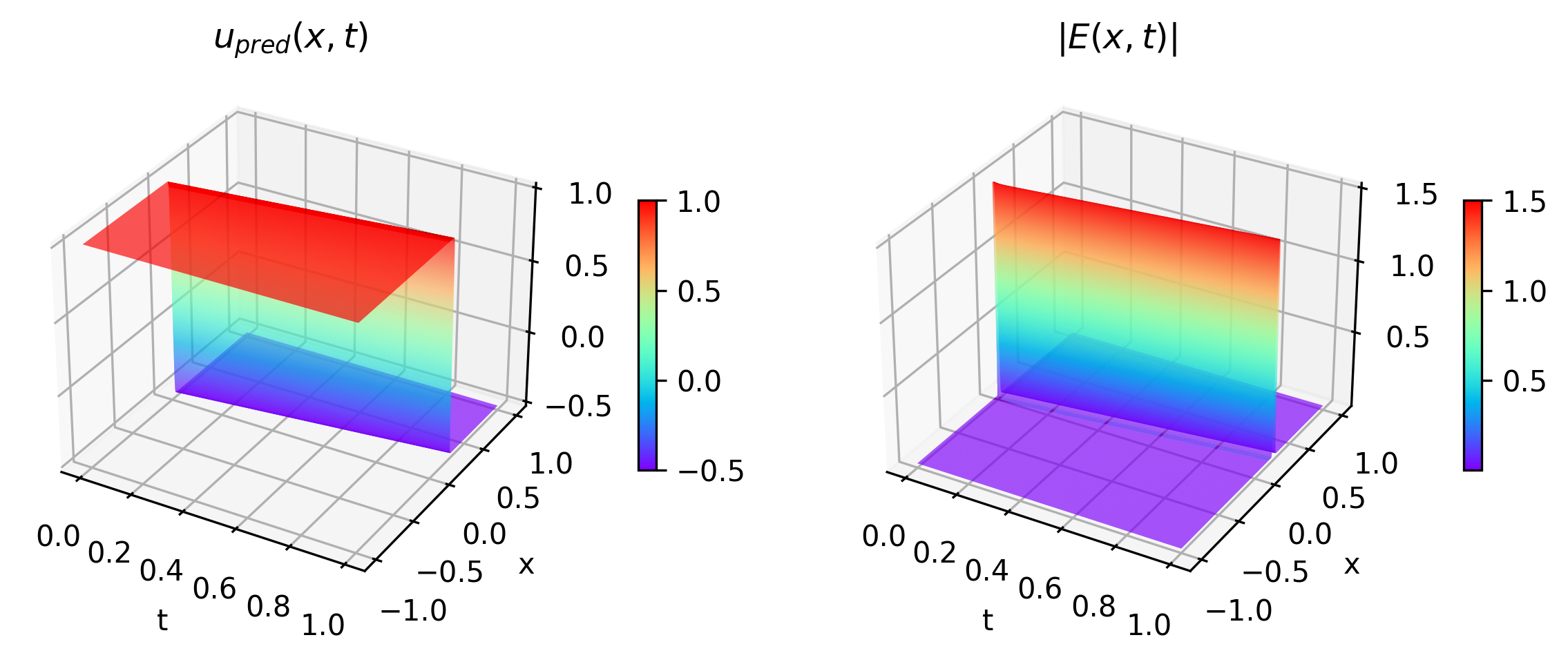}
        \caption{PINN prediction (left) and absolute pointwise error (right)}
    \end{subfigure}
    
    \bigskip
    
    \begin{subfigure}[b]{0.8\textwidth}
        \centering
        \includegraphics[width=\textwidth]{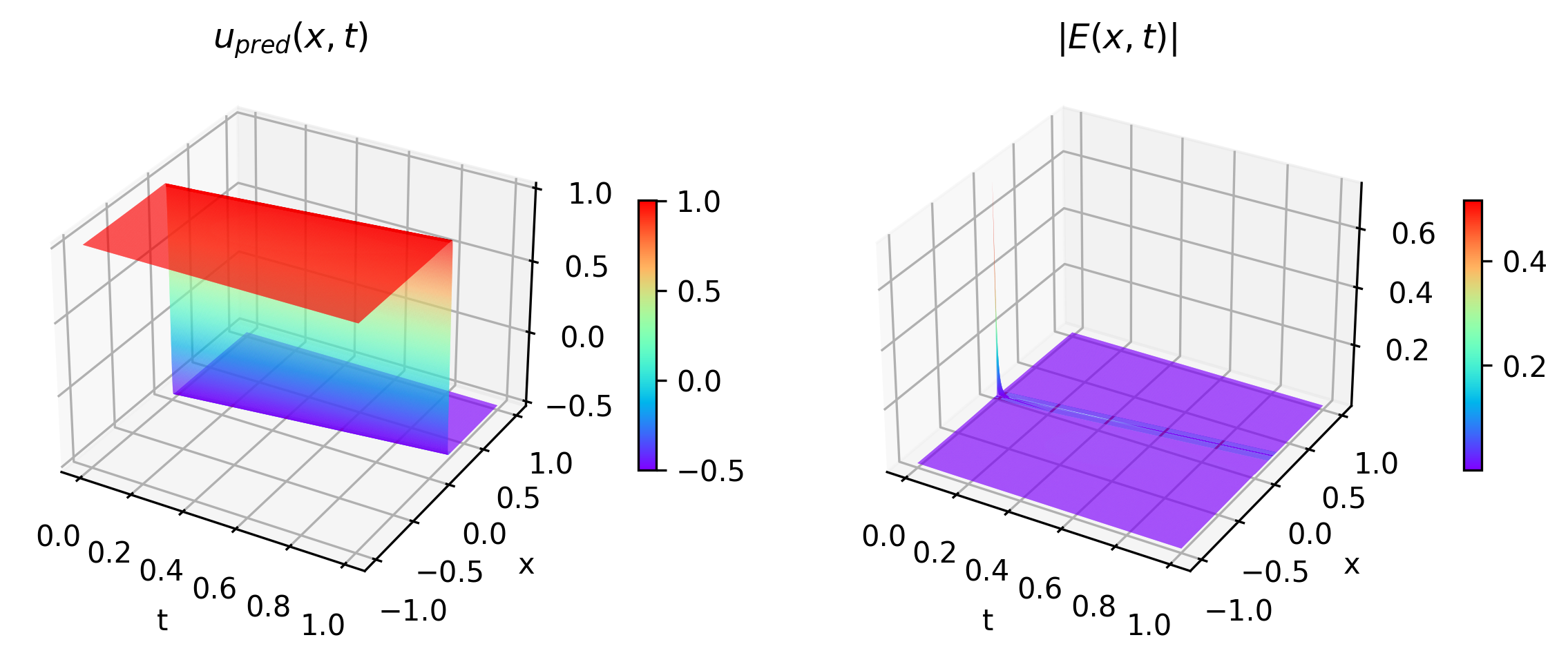}
        \caption{sl-PINN prediction (left) and absolute pointwise error (right)}
    \end{subfigure}
    \caption{Riemann problem moving shock case when $\ep=1/500$: PINN vs sl-PINN solution plots.}
    \label{plot2d nonsmooth 1}
\end{figure}

\begin{figure}
    \centering
    \begin{subfigure}[b]{0.4\textwidth}
        \centering
        \includegraphics[width=\textwidth]{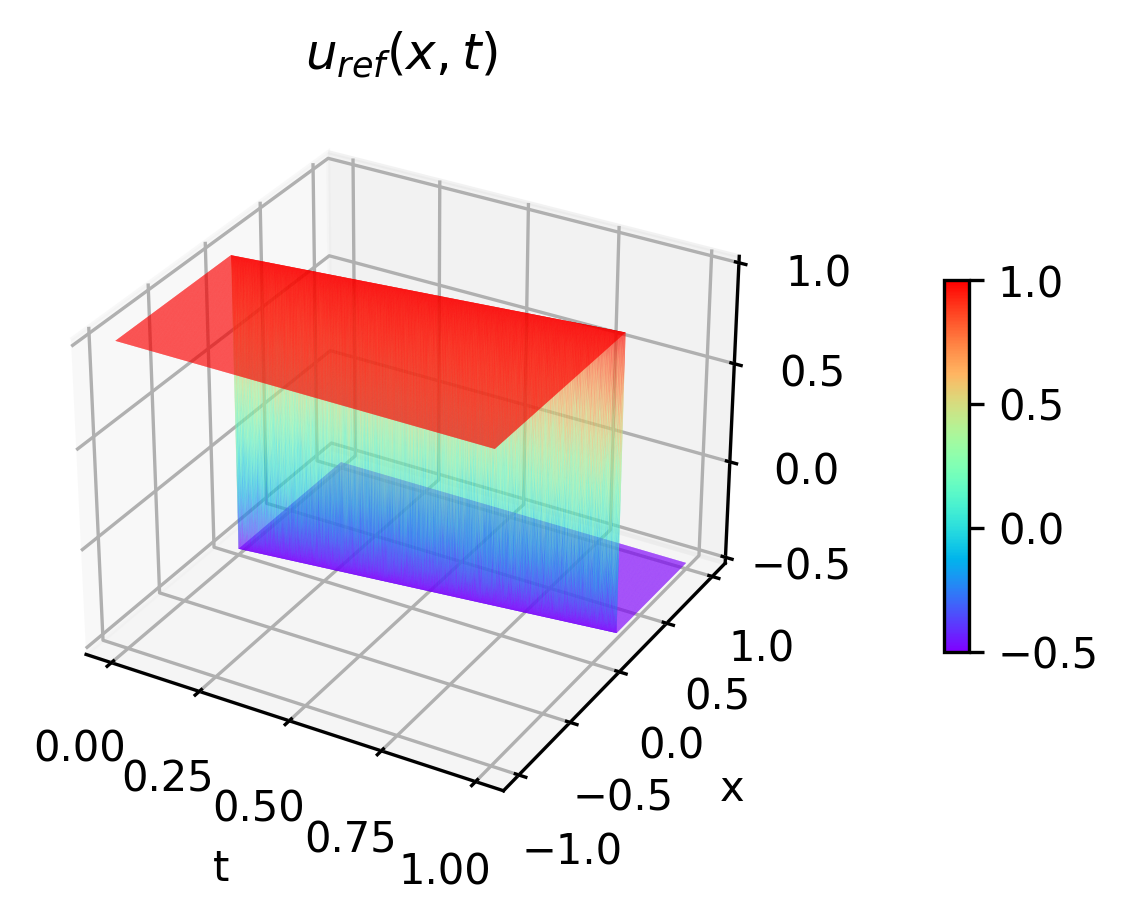}
        \caption{reference solution $\ep=1/10000$}
    \end{subfigure}
    
    \bigskip
    
    \begin{subfigure}[b]{0.8\textwidth}
        \centering
        \includegraphics[width=\textwidth]{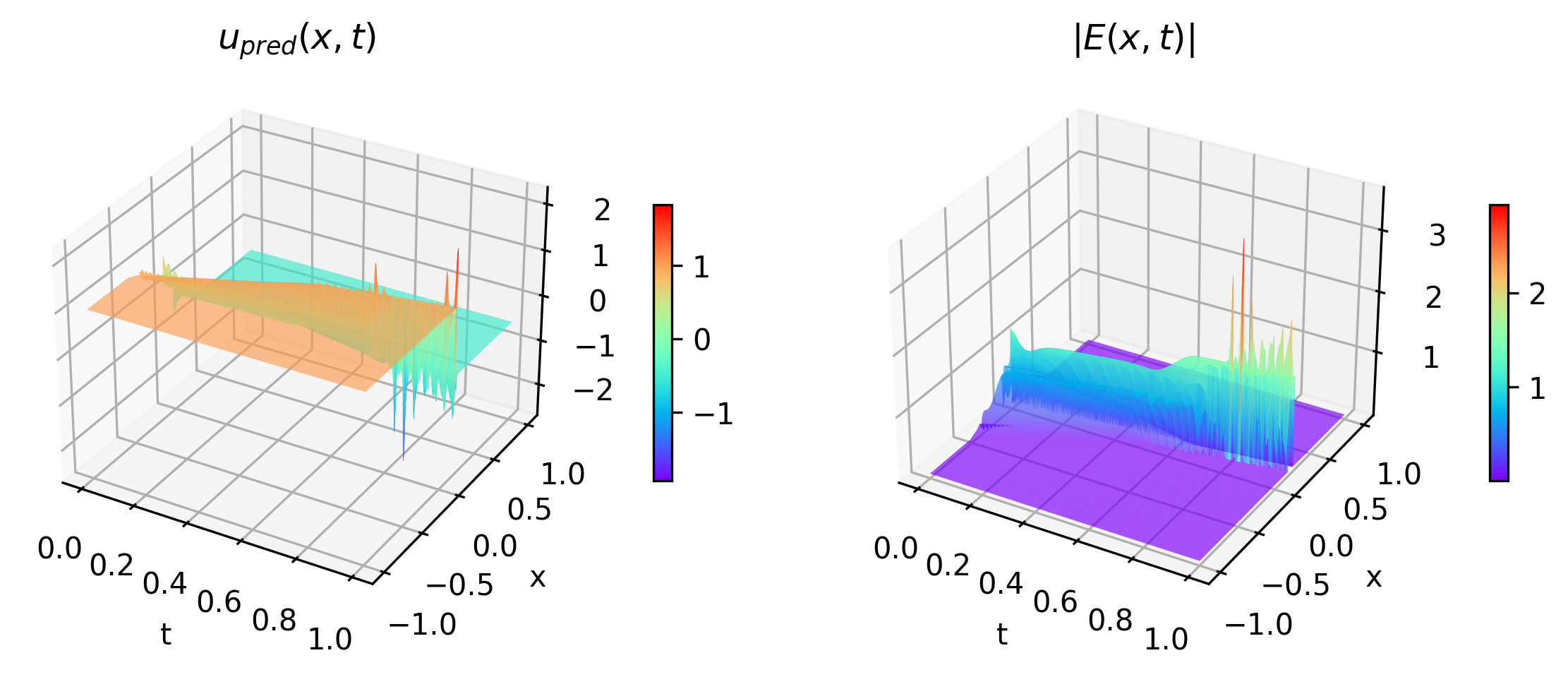}
        \caption{PINN prediction (left) and absolute pointwise error (right)}
    \end{subfigure}
    
    \bigskip
    
    \begin{subfigure}[b]{0.8\textwidth}
        \centering
        \includegraphics[width=\textwidth]{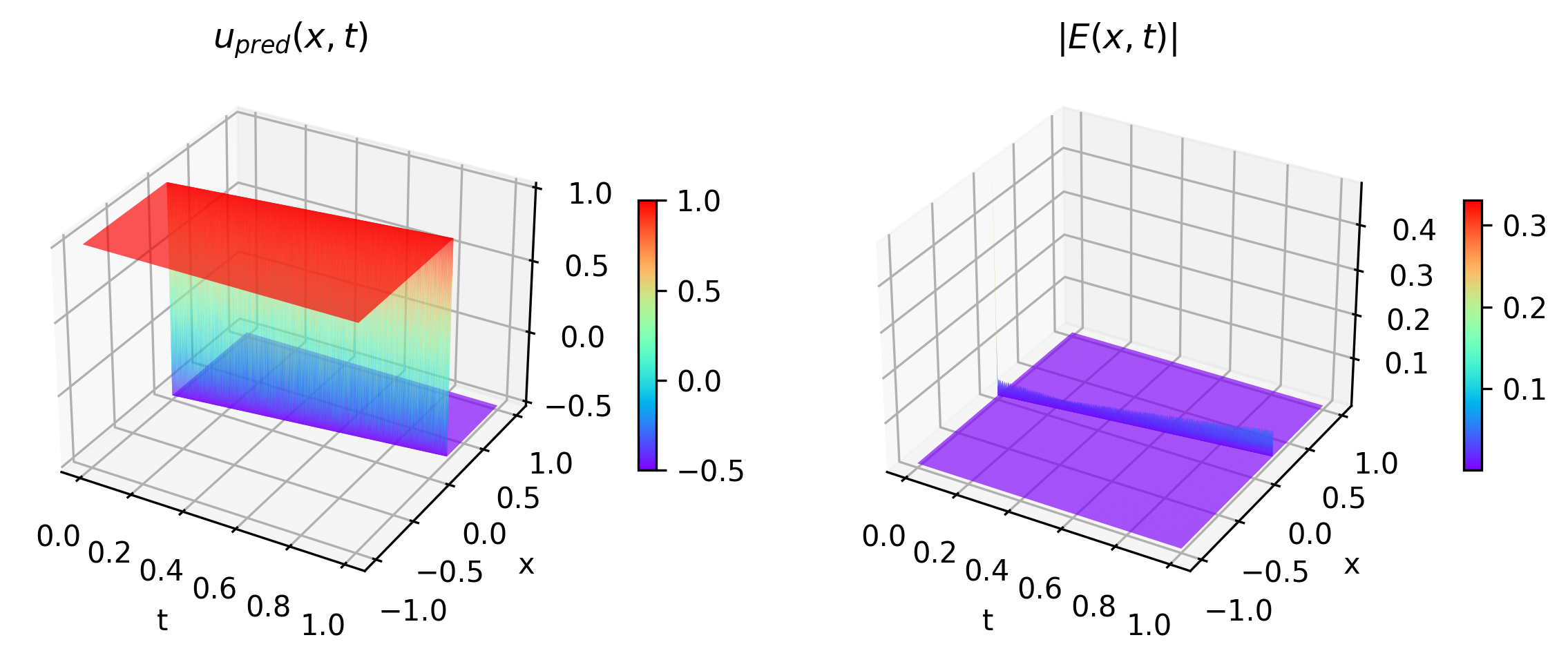}
        \caption{sl-PINN prediction (left) and absolute pointwise error (right)}
    \end{subfigure}
    \caption{Riemann problem moving shock case when $\ep=1/10000$: PINN vs sl-PINN solution plots.}
    \label{plot2d nonsmooth 2}
\end{figure}

\begin{figure}
    \centering
    \includegraphics[width=\textwidth]{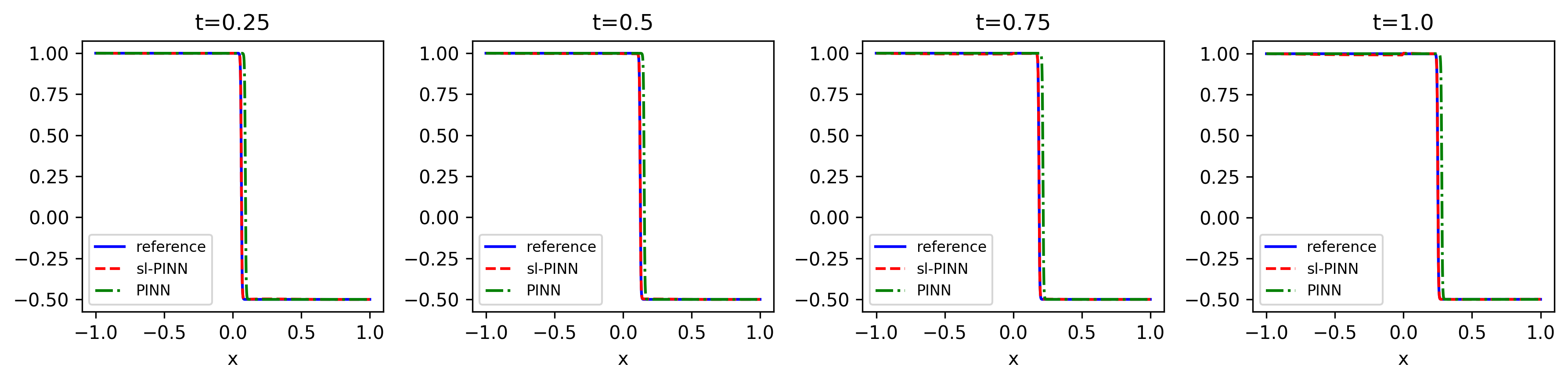}
    \caption{Riemann problem moving shock case: reference solution when $\ep=1/500$ vs PINN vs sl-PINN at specific time.}
    \label{plot1d nonsmooth combine 1}
\end{figure}

\begin{figure}
    \centering
    \begin{subfigure}[b]{\textwidth}
        \centering
        \includegraphics[width=\textwidth]{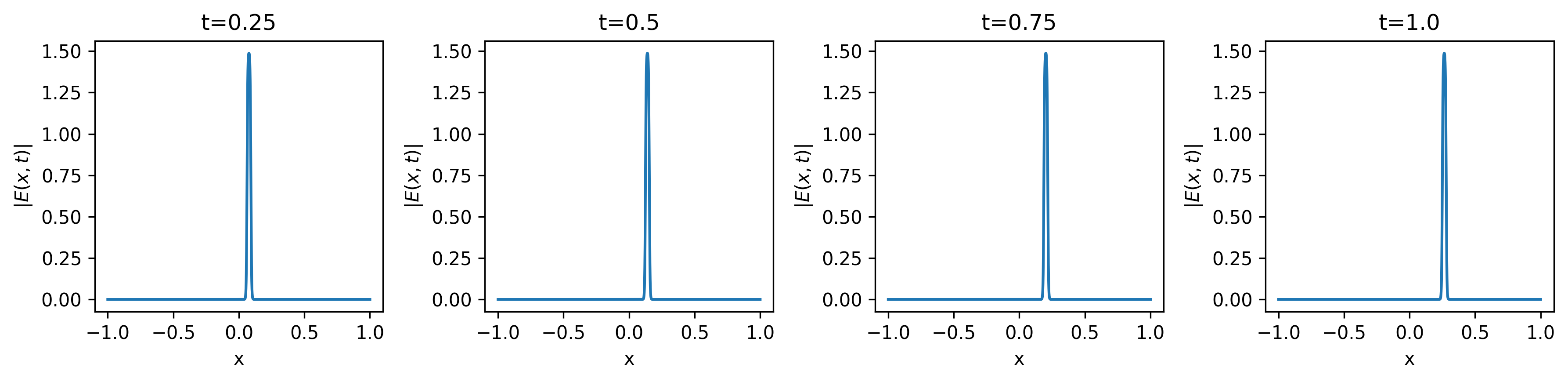}
        \caption{Errors of PINN prediction at specific time}
    \end{subfigure}
    
    \bigskip
    
    \begin{subfigure}[b]{\textwidth}
        \centering
        \includegraphics[width=\textwidth]{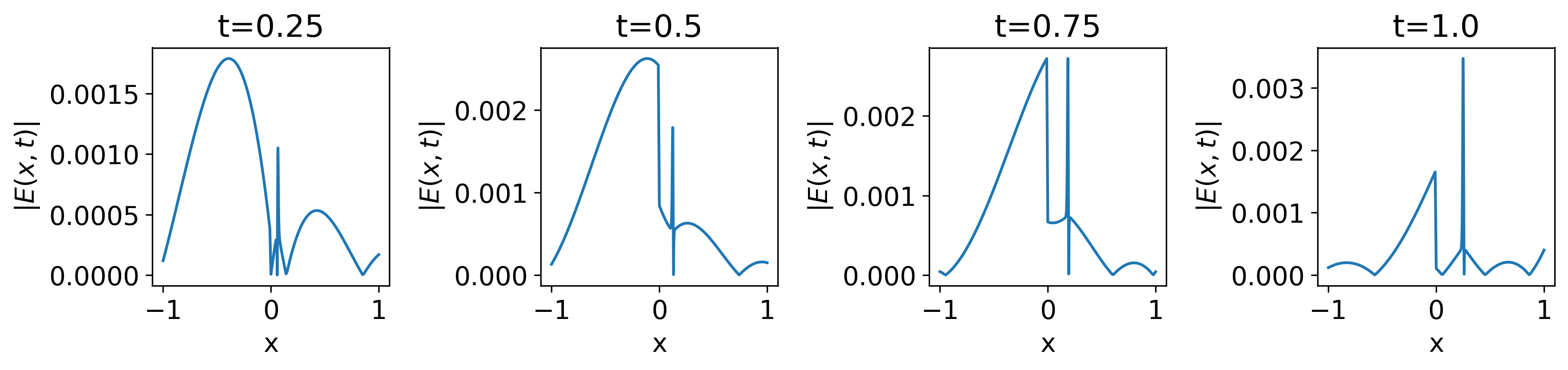}
        \caption{Errors of sl-PINN prediction at specific time}
    \end{subfigure}
    \caption{Riemann problem moving shock case when $\ep=1/500$: PINN vs sl-PINN solution error at specific time.}
    \label{plot1d nonsmooth 1}
\end{figure}

\begin{figure}
    \centering
    \includegraphics[width=\textwidth]{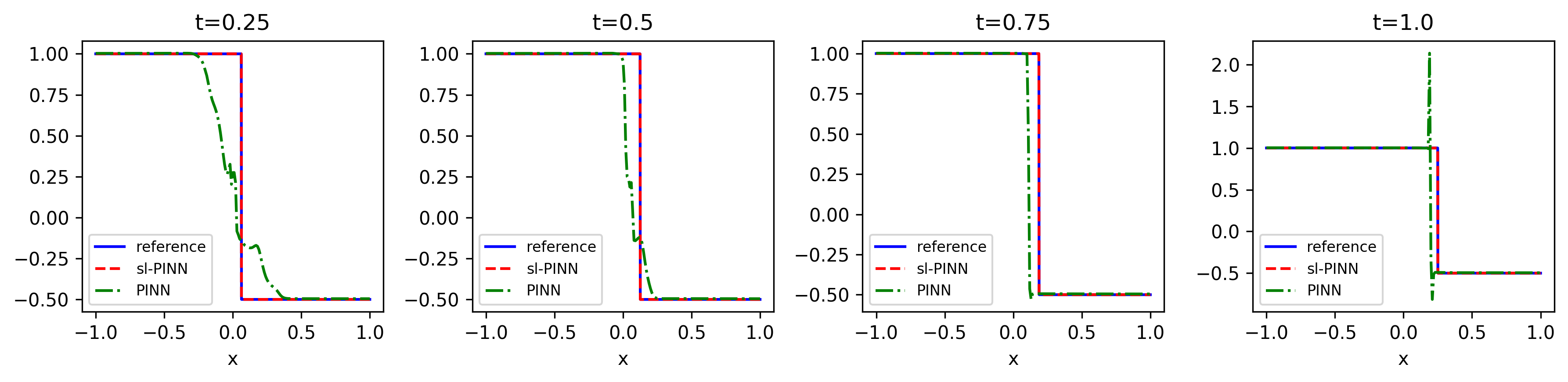}
    \caption{Riemann problem moving shock case when $\ep=1/10000$: reference solution vs PINN vs sl-PINN at specific time.}
    \label{plot1d nonsmooth combine 2}
\end{figure}

\begin{figure}
    \centering
    \begin{subfigure}[b]{\textwidth}
        \centering
        \includegraphics[width=\textwidth]{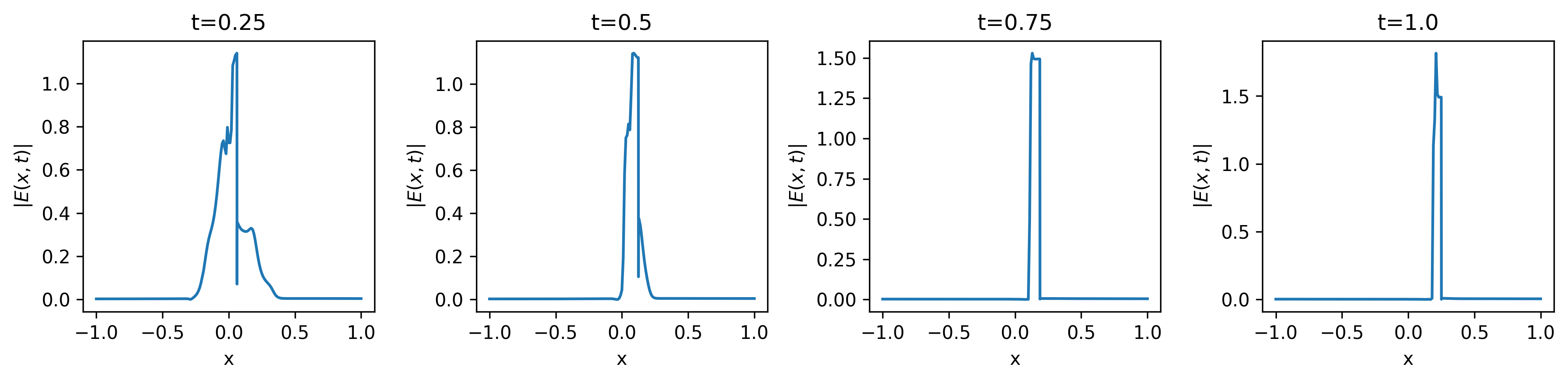}
        \caption{Errors of PINN prediction at specific time}
    \end{subfigure}
    
    \bigskip
    
    \begin{subfigure}[b]{\textwidth}
        \centering
        \includegraphics[width=\textwidth]{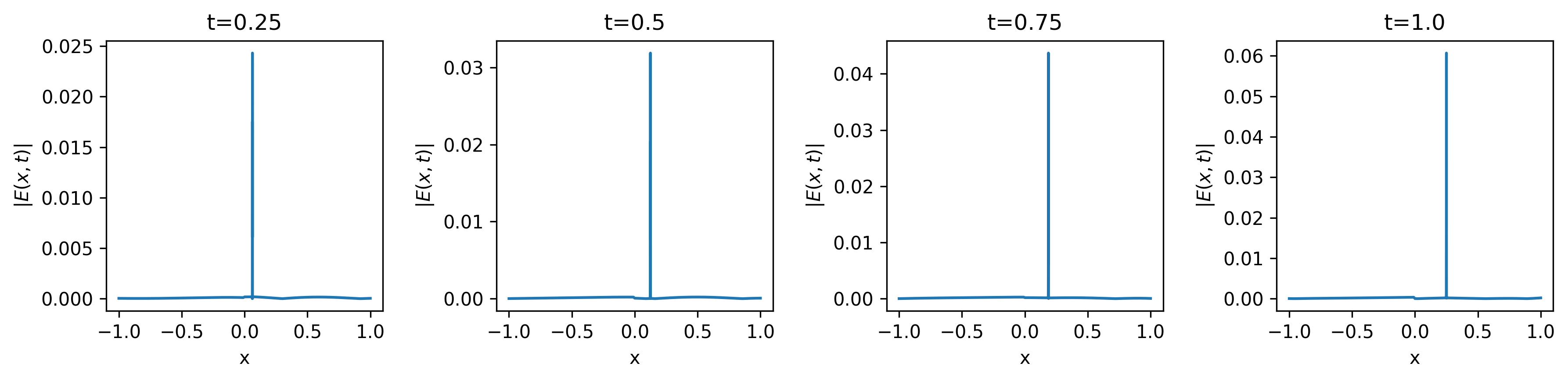}
        \caption{Errors of sl-PINN prediction at specific time}
    \end{subfigure}
    \caption{Riemann problem moving shock case when $\ep=1/10000$: PINN vs sl-PINN solution error at specific time.}
    \label{plot1d nonsmooth 2}
\end{figure}

After completing the training process, we investigate the accuracy in the $L^{2}$-norm (\ref{L2L2}) and the $L^{\infty}$-norm (\ref{Linf}) at specific time. 
The results of the steady shock case are shown in Table \ref{table riemann steady} and those of the moving shock case appear in Table \ref{table riemann moving}. 
We observe that for both cases, sl-PINN obtains accurate precision for all the small $\ep=1/500, \ 1/1000, \ 1/5000,$ and $1/10000$. On the contrary, the usual PINN has limited accuracy even for the case when a relatively large $\ep=1/500$ is considered. The solution plots of PINN and sl-PINN predictions for the moving shock case when $\ep=1/500$ and $\ep=1/10000$ appear in Figures \ref{plot2d nonsmooth 1} and \ref{plot2d nonsmooth 2}, their temporal snapshots in Figures \ref{plot1d nonsmooth combine 1} and \ref{plot1d nonsmooth combine 2}, and their error plots in Figures \ref{plot1d nonsmooth 1} and \ref{plot1d nonsmooth 2}. 
We find that sl-PINN is capable to capture the stiff behaviour of the Burgers' solution near the shock for a larger $\ep=1/500$ to a smaller $\ep=1/10000$.


We note that the initial data in these experiments are discontinuous. 
It is well-known that neural networks (or any traditional numerical method) have limited performance in fitting a target function with discontinuity \cite{JKK20, RB19}. Therefore, as seen in the predictions of the sl-PINN in Figure \ref{plot2d nonsmooth 1} or \ref{plot2d nonsmooth 2}, a large error near the singular layer is somewhat expected.

\begin{table}[]
\centering
\begin{subtable}{\textwidth}
\centering
\begin{tabular}{|l|r|r|}
\hline
$\| E \|_{L^{2}_{t}( L^{2}_{x}(\Omega))}$ & \multicolumn{1}{l|}{PINNs} & \multicolumn{1}{l|}{sl-PINNs} \\ \hline
$\ep=1/500$                               & 6.77E-01                   & 4.34E-03                      \\
$\ep=1/1000$                              & 1.08E+00                   & 4.28E-03                      \\
$\ep=1/5000$                              & 6.84E-01                   & 2.70E-03                      \\
$\ep=1/10000$                             & 9.29E-01                   & 3.74E-03                      \\ \hline
\end{tabular}
\bigskip
\end{subtable}
\begin{subtable}{\textwidth}
\centering
\begin{tabular}{|l|rr|rr|}
\hline
                                   & \multicolumn{2}{c|}{PINNs}                                            & \multicolumn{2}{c|}{sl-PINNs}                                         \\ \hline
$\| E \|_{L^{\infty}_{x}(\Omega)}$ & \multicolumn{1}{l|}{$\ep=1/500$} & \multicolumn{1}{l|}{$\ep=1/10000$} & \multicolumn{1}{l|}{$\ep=1/500$} & \multicolumn{1}{l|}{$\ep=1/10000$} \\ \hline
$t = 0.25$                         & \multicolumn{1}{r|}{2.00E+00}    & 2.00E+00                           & \multicolumn{1}{r|}{3.82E-04}    & 4.80E-02                           \\
$t = 0.5$                          & \multicolumn{1}{r|}{2.00E+00}    & 1.99E+00                           & \multicolumn{1}{r|}{4.64E-04}    & 4.79E-02                           \\
$t = 0.75$                         & \multicolumn{1}{r|}{2.00E+00}    & 1.99E+00                           & \multicolumn{1}{r|}{4.12E-04}    & 4.80E-02                           \\
$t = 1.0$                          & \multicolumn{1}{r|}{2.00E+00}    & 1.90E+00                           & \multicolumn{1}{r|}{6.11E-04}    & 4.80E-02                           \\ \hline
\end{tabular}
\end{subtable}
\caption{Riemann problem - steady shock case $a=1$, $b=1$: $L^{2}$-error comparison (top) and $L^{\infty}$-error comparison (bottom)}
\label{table riemann steady}
\end{table}

\begin{table}[]
\centering
\begin{subtable}{\textwidth}
\centering
\begin{tabular}{|l|r|r|}
\hline
$\| E \|_{L^{2}_{t}( L^{2}_{x}(\Omega))}$ & \multicolumn{1}{l|}{PINNs} & \multicolumn{1}{l|}{sl-PINNs} \\ \hline
$\ep=1/500$                               & 5.84E-01                   & 5.46E-03                      \\
$\ep=1/1000$                              & 6.23E-01                   & 4.19E-03                      \\
$\ep=1/5000$                              & 4.24E-01                   & 3.01E-03                      \\
$\ep=1/10000$                             & 4.04E-01                   & 3.84E-03                      \\ \hline
\end{tabular}
\bigskip
\end{subtable}
\begin{subtable}{\textwidth}
\centering
\begin{tabular}{|l|rr|rr|}
\hline
                                   & \multicolumn{2}{c|}{PINNs}                                            & \multicolumn{2}{c|}{sl-PINNs}                                         \\ \hline
$\| E \|_{L^{\infty}_{x}(\Omega)}$ & \multicolumn{1}{l|}{$\ep=1/500$} & \multicolumn{1}{l|}{$\ep=1/10000$} & \multicolumn{1}{l|}{$\ep=1/500$} & \multicolumn{1}{l|}{$\ep=1/10000$} \\ \hline
$t = 0.25$                         & \multicolumn{1}{r|}{1.49E+00}    & 1.14E+00                           & \multicolumn{1}{r|}{1.79E-03}    & 2.43E-02                           \\
$t = 0.5$                          & \multicolumn{1}{r|}{1.49E+00}    & 1.14E+00                           & \multicolumn{1}{r|}{2.63E-03}    & 3.19E-02                           \\
$t = 0.75$                         & \multicolumn{1}{r|}{1.49E+00}    & 1.53E+00                           & \multicolumn{1}{r|}{2.72E-03}    & 4.37E-02                           \\
$t = 1.0$                          & \multicolumn{1}{r|}{1.49E+00}    & 1.81E+00                           & \multicolumn{1}{r|}{3.47E-03}    & 6.07E-02                           \\ \hline
\end{tabular}
\end{subtable}
\caption{Riemann problem - moving shock case $a=1$, $b=0.5$: $L^{2}$-error comparison (top) and $L^{\infty}$-error comparison (bottom)}
\label{table riemann moving}
\end{table}

\section{Conclusion}\label{S.conclusion}

In this article, we introduce a new learning method known as sl-PINN to address the one-dimensional viscous Burgers problem with low viscosity, causing a singular interior layer. The method incorporates corrector functions that describe the robust behavior of the solution near the interior layer into the structure of PINNs to improve learning performance in that area. This is achieved through interior analysis, treating the problem as two boundary layer sub-problems. The profile of the corrector functions is studied using asymptotic analysis. Both stationary and moving shock cases have been explored. Numerical experiments show that our sl-PINNs accurately predict the solution for every small viscosity, which, in particular, reduces errors near the interior layer compared to the original PINNs. Our proposed method provides a comprehensive understanding of the solution’s behavior near the interior layer, aiding in capturing the robust part of the training solution. This approach offers a distinct perspective from other existing works \cite{SWYTC24, ZYZK24, LMMK21}, achieving better accuracy, especially when the shock is more robust.


\section*{Acknowledgment}
T.-Y. Chang gratefully was supported by
the Graduate Students Study Abroad Program grant funded by the National Science and Technology Council(NSTC) in Taiwan. 
G.-M. Gie was supported by 
Simons Foundation Collaboration Grant for Mathematicians. 
Hong was supported by 
Basic Science Research Program through the National Research Foundation of Korea (NRF) funded by the Ministry of Education (NRF-2021R1A2C1093579) and by the Korea government(MSIT) (RS-2023-00219980).  
Jung was supported by 
National Research Foundation of Korea(NRF) grant funded by the Korea government(MSIT) (No. 2023R1A2C1003120).

\end{document}